\documentclass[final,onefignum,onetabnum]{siamart190516}

\ifpdf
\hypersetup{
  pdftitle={AL Preconditioner for Implicitly Constituted Flow},
  pdfauthor={P.A. Gazca--Orozco, P.E. Farrell}
}
\fi

\usepackage{amsmath,amssymb,amsbsy,amsfonts}
\usepackage{mathtools}
\usepackage{bm}
\usepackage{stmaryrd}
\usepackage{booktabs}
\usepackage{multirow}
\usepackage{graphicx}
\usepackage{epstopdf}
\ifpdf
  \DeclareGraphicsExtensions{.eps,.pdf,.png,.jpg}
\else
  \DeclareGraphicsExtensions{.eps}
\fi
\usepackage{tikz}
\usepackage{pgfplots}
\pgfplotsset{width=7cm,compat=1.14}
\usepackage[caption=false]{subfig}
\def\abs#1{\left| #1 \right|}

\DeclarePairedDelimiter\jump{\llbracket}{\rrbracket}
\def\tens#1{\pmb{\mathsf{#1}}}
\def\vec#1{\boldsymbol{#1}}
\def\R{\mathbb{R}}

\def\sym{\mathop{\mathrm{sym}}\nolimits}
\def\tr{\mathop{\mathrm{tr}}\nolimits}
\def\Rds{\mathbb{R}^{d \times d}_{\sym}}
\def\Rdst{\mathbb{R}^{d \times d}_{\sym,\tr}}
\def\supp{\mathop{\mathrm{supp}}\nolimits}

\def\lin{\mathop{\mathrm{span}}\nolimits} 
\def\macrostar{\mathop{\mathrm{macrostar}}\nolimits} 
\def\diver{\mathop{\mathrm{div}}\nolimits} 

\def\Du{\BD(\bu)}

\def\Lsymtr{L_{\sym,\tr}}
\def\e{\mathrm{e}}
\def\b0{\vec{0}}

\def\bf{\vec{f}}

\def\bu{\vec{u}}
\def\bv{\vec{v}}
\def\bw{\vec{w}}

\def\bz{\vec{z}}

\def\bphi{\vec{\varphi}}

\def\bsigma{\vec{\sigma}}
\def\btau{\vec{\tau}}
\def\B0{\tens{0}}

\def\BD{\tens{D}}

\def\BG{\tens{G}}

\def\BI{\tens{I}}

\def\BS{\tens{S}}

\newsiamremark{remark}{Remark}
\newsiamremark{hypothesis}{Hypothesis}
\newsiamremark{example}{Example}
\crefname{hypothesis}{Hypothesis}{Hypotheses}
\newsiamthm{claim}{Claim}
\newsiamthm{assumption}{Assumption}
\crefname{assumption}{Assumption}{Assumptions}

\headers{AL Preconditioner for Implicitly Constituted Flow}{P. E. Farrell, and  P. A. Gazca-Orozco}

\title{An augmented Lagrangian preconditioner for implicitly-constituted non-Newtonian incompressible flow \thanks{Submitted to the editors May 2020.
\funding{This research is supported by the Engineering and Physical Sciences
Research Council [grant numbers EP/R029423/1 and EP/V001493/1], and by the
EPSRC Centre For Doctoral Training in Partial Differential Equations:
Analysis and Applications [grant number EP/L015811/1]. The second author was supported by CONACyT (Scholarship 438269).}}}

\author{P.\ E.\ Farrell\thanks{Mathematical Institute, University of Oxford, UK (\email {patrick.farrell@maths.ox.ac.uk}).}
\and P.\ A.\ Gazca-Orozco\thanks{Mathematical Institute, University of Oxford, UK (\email {gazcaorozco@maths.ox.ac.uk}).}
}

\usepackage{amsopn}

\begin{document}

\maketitle

\begin{abstract}
We propose an augmented Lagrangian preconditioner for a three-field stress-velocity-pressure discretization of stationary non-Newtonian incompressible flow with an implicit constitutive relation of power-law type.
	The discretization employed makes use of the divergence-free Scott--Vogelius pair for the velocity and pressure. The preconditioner builds on the work [P.~E.~Farrell, L.~Mitchell, and F.~Wechsung, \textit{SIAM J.\ Sci.\ Comput.}, 41 (2019), pp. A3073--A3096], where a Reynolds-robust preconditioner for the three-dimensional Newtonian system was introduced. The preconditioner employs a specialized multigrid method for the stress-velocity block that involves a divergence-capturing space decomposition and a custom prolongation operator. The solver exhibits excellent robustness with respect to the parameters arising in the constitutive relation, allowing for the simulation of a wide range of materials.
\end{abstract}



\section{Introduction}

For $d\in\{2,3\}$ let $\Omega\subset \mathbb{R}^d$ be a bounded polygonal domain with Lipschitz boundary. The goal of this work is to construct a preconditioner for the Newton linearization of a system describing the steady state of an incompressible fluid:
\begin{subequations}\label{eq:PDE_continuous}
	\begin{alignat}{2}
		\begin{aligned}
			-\diver \BS + \diver (\bu\otimes&\bu) + \nabla p = \bm{f} \quad & &\text{ in }\Omega,\\
			\diver\bu &= 0\quad & &\text{ in }\Omega,\\
			\bu &= \bu_0\quad & &\text{ on }\partial\Omega,\\
		\end{aligned}
	\end{alignat}
	where $\bm{f}$ and $\bu_0$ are given. In the equation above, $\bu\colon \Omega\to \R^d$ denotes the velocity field, $p:\Omega\to \R$ denotes the pressure (mean normal stress) and $\BS:\Omega\to \Rdst$ is the shear stress. In order to ensure the uniqueness of the pressure we impose a zero mean constraint $\int_\Omega p = 0$. The system is closed with an implicit constitutive relation of the form 
\begin{alignat}{2}\label{eq:ImplicitCR}
	\begin{aligned}
	\BG(\cdot,\BS,\BD):=\alpha(\cdot,|\BS|^2,|\BD|^2)\BD - \beta(\cdot,|\BS|^2,|\BD|^2)\BS = \b0,
	\end{aligned}
\end{alignat}
\end{subequations}
where $\BD := \Du = \frac{1}{2}(\nabla\bu+\nabla\bu^{\top})$ is the symmetric velocity gradient, $\alpha,\beta:\Omega\times[0,\infty)^2\rightarrow \mathbb{R}^+$ are positive functions,  and $\BG:\Omega \times \Rds\times \Rds \rightarrow \Rds$ is a function that defines a monotone graph; this includes for instance the usual Navier--Stokes and power-law models (the precise assumptions will be introduced later). This framework of implicitly constituted fluids is very natural when modelling wide classes of materials, and allows for their systematic study in a thermodynamically consistent manner \cite{Rajagopal:2003,Rajagopal2006,Rajagopal2008}.
Since in general it is not possible to solve explicitly for the stress $\BS$ in \eqref{eq:ImplicitCR} to substitute it into the momentum equation, we consider a three-field formulation of the problem in which the stress is one of the unknowns. After discretization and Newton linearization, the system has the following block form:
\begin{equation}\label{eq:block2}
    \begin{bmatrix}
    A & B^\top\\ B & 0
    \end{bmatrix}
    \begin{bmatrix}
    \bz \\ p
    \end{bmatrix}
    =
    \begin{bmatrix}
    \bf \\ g
    \end{bmatrix},
\end{equation}
%
where $\bz := (\BS,\bu)^\top$, $A$ is the stress-velocity block and $B$ represents the discrete divergence on the velocity space (c.f.\ \eqref{eq:block_structure} below). A popular approach to preconditioning systems with this structure is based on the block factorization
\begin{equation*}
\begin{bmatrix}
A & B^{\top} \\B &0
\end{bmatrix}^{-1} =  \begin{bmatrix}
I & -A^{-1}B^{\top}\\ 0 & I
\end{bmatrix} \begin{bmatrix}
A^{-1} & 0 \\ 0 & S^{-1}
\end{bmatrix}\begin{bmatrix}
I & 0 \\ -BA^{-1} & I
\end{bmatrix},
\end{equation*}
where $S = -BA^{-1} B^\top$ is the Schur complement. If approximations $\tilde{A}^{-1}$ and $\tilde{S}^{-1}$ of $A^{-1}$ and $S^{-1}$ are available, they can be used in this formula to precondition the coupled system. For a velocity-pressure formulation of the Stokes system, it is known that the Schur complement is spectrally equivalent to the viscosity-weighted pressure mass matrix \cite{Silvester:1994,Mardal2011}: $S \sim -\nu^{-1} M_p$. Using (for instance) an algebraic multigrid cycle on $A$ as $\tilde{A}^{-1}$ and the inverse diagonal of the pressure mass matrix as $\tilde{S}^{-1}$ results in a mesh-independent preconditioner for the Stokes system. For the Navier--Stokes system this choice results in a solver whose performance degrades badly as the Reynolds number $\mathrm{Re}$ increases, i.e.\ the number of Krylov iterations per nonlinear iteration grows with $\mathrm{Re}$ \cite{Elman1996}. Other preconditioners such as the pressure convection-diffusion \cite{Kay:2006} and least-squares commutator \cite{Elman2006} perform well for moderate Reynolds numbers, but their performance still deteriorates as the Reynolds number grows \cite{Elman2014}.

An alternative approach for dealing with the Schur complement approximation was proposed by Benzi and Olshanskii \cite{Benzi2006} for a 2D Navier--Stokes problem and later extended to the 3D problem by Farrell, Mitchell and Wechsung \cite{Farrell2019a}. The main idea is to modify the system by adding an augmented Lagrangian term:
\begin{equation}\label{eq:AL_stabilisation}
\begin{bmatrix}
A + \gamma B^{\top}M_p^{-1}B & B^{\top}\\B &0
\end{bmatrix}
\begin{bmatrix}
\bz \\ p
\end{bmatrix}
=
\begin{bmatrix}
\bf + \gamma B^{\top}M_p^{-1} g \\ g
\end{bmatrix},
\end{equation}
where $\gamma>0$ is a parameter. Observe that this modification does not change the solution of the system, since $B\bz=g$. The continuous form of the term $\gamma B^{\top}M_p^{-1}B$ could be interpreted as a least-squares term that penalizes the $L^2$ norm of $\diver\bu$, and appears in other contexts, such as the iterated penalty and artificial compressibility methods \cite{Temam1968,Chorin1967}. From the Sherman--Morrison--Woodbury formula (see e.g.\ \cite{Bacuta2006}) we see that the inverse Schur complement of the augmented matrix can be approximated as 
\begin{align*}
	S^{-1} &= (-B(A + \gamma B^{\top}M_p^{-1}B)^{-1}B^{\top})^{-1} = -(BA^{-1}B^{\top})^{-1} - \gamma M_p^{-1} \\
	       &\approx -(\nu + \gamma)M_p^{-1} \approx -\gamma M_p^{-1},
\end{align*}
with the approximation improving as $\gamma \rightarrow \infty$ (cf.\ \cite{Farrell2019a}).

The challenge is to develop an efficient solver for the augmented \textcolor{black}{$(1,1)$ block} $A+ \gamma B^{\top}M_p^{-1}B$. This is not trivial as the augmented Lagrangian term has a large kernel (all divergence-free velocity fields) and so the matrix degenerates as $\gamma \rightarrow \infty$. The essential breakthrough for the Navier--Stokes system came with the work \cite{Benzi2006}, where a specialized multigrid operator was developed for the \textcolor{black}{$(1,1)$ block}, applying ideas developed by Sch\"{o}berl for nearly incompressible elasticity \cite{Schoeberl:1999,Schoeberl1999}. In this work we will apply these ideas to develop a robust multigrid operator for the coupled stress-velocity block. The two main components needed to obtain a robust multigrid solver are a robust relaxation and a robust prolongation operator, which we will develop in the following sections. \textcolor{black}{In previous work Farrell et al.~developed a preconditioner in this framework for the Scott--Vogelius discretization of the Newtonian case \cite{Farrell2020}, with a specialised multigrid method applied to the velocity problem arising after augmentation; the main challenge in this work is the development of an appropriate inner solver for the augmented stress-velocity block that is required in the implicitly constituted non-Newtonian case. This inner system presents a saddle point structure of its own, which we tackle with suitable monolithic multigrid techniques.}

It is important to note that the available theory for the development of robust multigrid solvers assumes that the matrix $A$ is symmetric and positive definite (SPD). This assumption does not hold for the problem under consideration; the stress-velocity block in \eqref{eq:block2} itself has a saddle point structure and is not symmetric. Nevertheless, satisfying the requirements of the SPD case appears to give good performance in the
general case also, as observed in the computational experiments of previous works \cite{Benzi2006,Farrell2019a,Farrell2020}. The computational experiments of \Cref{Sec:Examples} demonstrate that the preconditioner we propose possesses similarly excellent robustness with respect to parameters arising in the implicit constitutive relation \eqref{eq:ImplicitCR}.

\subsection{Implicit constitutive relation}
We will employ standard notation for Lebesgue and Sobolev spaces (e.g.\ $(W^{k,r}(\Omega),\|\cdot\|_{W^{k,r}(\Omega)})$ and $(L^q(\Omega),\|\cdot\|_{L^q(\Omega)})$) \textcolor{black}{(see e.g. \cite{Adams2003})}. The space $W^{k,r}_0(\Omega)$ is defined for $r\in [1,\infty)$ as the closure of the space of smooth functions with compact support $C_0^\infty(\Omega)$ with respect to the norm $\|\cdot\|_{W^{k,r}(\Omega)}$ and  we will denote the dual space of $W^{1,r}_0(\Omega)$ by $W^{-1,r'}(\Omega)$. Here $r'$ denotes the H\"{o}lder conjugate of $r$, i.e.\ the number defined by $1/r + 1/r' =1$. The space of traces on the boundary of functions in $W^{1,r}(\Omega)$ will be denoted by $W^{1/r',r}(\partial\Omega)$.  For $r\in [1,\infty)$ we also define the following useful subspaces:
\begin{gather*}
L^r_0(\Omega) := \left\{q\in L^r(\Omega) \, : \, \int_\Omega q = 0\right\},\\
W^{1,r}_{0,{{\diver}}}(\Omega)^d := \overline{\{\bv\in C^\infty_0(\Omega)^d\, :\, {{\diver}}\,\bv=0\}}^{\|\cdot\|_{W^{1,r}(\Omega)}},\\
L_{{{\tr}}}^r(Q)^{d\times d} := \{\bm{\tau}\in L^r(Q)^{d\times d}\, : \, {{\tr}}(\bm{\tau}) = 0\},\\
L_{{\sym}}^r(Q)^{d\times d} := \{\bm{\tau}\in L^r(Q)^{d\times d}\, : \, \bm{\tau}^\top = \bm{\tau}\},\\
L_{{\sym,\tr}}^r(Q)^{d\times d} := L_{{\sym}}^r(Q)^{d\times d} \cap L_{{{\tr}}}^r(Q)^{d\times d}.
\end{gather*}

In the definition of the space $L_{{{\tr}}}^r(Q)^{d\times d}$ above, ${{\tr}}(\bm{\tau})$ denotes the usual matrix trace of the $d\times d$ matrix function $\bm{\tau}$.

We will assume that the function $\BG$ in \eqref{eq:ImplicitCR} satisfies the following conditions:

\begin{enumerate}
	\item[(A1)] The mapping $(\BD,\BS)\in \Rds\times\Rds \mapsto \BG(x,\BS,\BD)$ is Fr\'{e}chet-differentiable for almost every $x\in\Omega$.
	\item[(A2)] The mapping $x\in\Omega \mapsto \BG(x,\BS,\BD)$ belongs to $L^\infty(\Omega;\Rds)$ for every $(\BD,\BS)\in \Rds\times \Rds$.
\end{enumerate}

The differentiability assumption (A1) is needed because Newton's method will be applied to linearize the system. Let us introduce the graph $\mathcal{A}:\Omega \rightarrow \Rds\times\Rds$ defined by $\BG$ using the canonical identification 
\begin{equation}
(\BD, \BS)\in\mathcal{A}(\cdot) \Longleftrightarrow \BG(\cdot,\BS,\BD)=\b0.
\end{equation}

We will assume in addition that, for some $r >1$, $\mathcal{A}(\cdot)$ satisfies the following properties for almost every $x\in\Omega$:
\begin{enumerate}
\item[(A3)] [$\mathcal{A}$\textit{ is a strictly monotone graph}] For every $(\BD_1,\BS_1), (\BD_2,\BS_2)\in \mathcal{A}(x)$,
\[(\BS_1 - \BS_2):(\BD_1 -\BD_2)\geq 0,\]
with strict inequality when $(\BD_1,\BS_1)\neq(\BD_2,\BS_2)$.
\item[(A4)] [$\mathcal{A}$\textit{ is an }$r$\textit{-graph}] There is a non-negative function $m\in L ^1(\Omega)$ and a constant $c>0$ such that
\[\BS:\BD \geq -m + c(|\BD|^r + |\BS|^{r'})\quad \text{for all } (\BD,\BS)\in\mathcal{A}(x).\]
\end{enumerate}
For a constitutive relation of the form \eqref{eq:ImplicitCR}, the assumptions above have the following important consequences, for almost every $x\in\Omega$:
\begin{enumerate}
\item[(C1)] [$\mathcal{A}$\textit{ includes the origin}] $(\mathbf{0},\mathbf{0})\in \mathcal{A}(x)$.
\item[(C2)] [$\mathcal{A}$\textit{ is maximal monotone}] If $(\BD,\BS) \in \Rds\times \Rds$ is such that
\[(\hat{\BS} - \BS):(\hat{\BD} -\BD)\geq 0 \quad \text{for all }(\hat{\BD},\hat{\BS})\in \mathcal{A}(x),\]
then $(\BD,\BS) \in \mathcal{A}(x)$.
\item[(C3)] [\textit{Measurability}] \textcolor{black}{The set-valued map $x\mapsto \mathcal{A}(x)$ is $\mathcal{L}(\Omega)$--$(\mathcal{B}(\Rds)\otimes\mathcal{B}( \Rds))$ measurable; here $\mathcal{L}(\Omega)$ denotes the family of Lebesgue measurable subsets of $\Omega$ and $\mathcal{B}(\Rds)$ is the family of Borel subsets of $\Rds$.}
\item[(C4)] [\textit{Compatibility}] For any $(\BD,\BS)\in\mathcal{A}(x)$ we have that
\begin{equation*}
\tr(\BD) = 0 \Longleftrightarrow \tr(\BS)=0.
\end{equation*}
\end{enumerate}
\noindent
In other words, we have that $\mathcal{A}$ is a \emph{maximal monotone $r$-graph}. In the weak formulation of this system we look for $(\BS,\bu,p)\in \Lsymtr^{r'}(\Omega)^{d\times d}\times\left(\bu_0 + W^{1,r}_{0,\diver}(\Omega)^d\right)\times L^{\tilde{r}}_0(\Omega)$, where $\tilde{r}:= \min\{r',r^*/2\}$, such that 

\begin{subequations}\label{eq:WeakFormulation}
\begin{alignat}{2}
\int_\Omega \BS:\BD(\bv) -\int_\Omega \bu\otimes\bu&:\BD(\bv) -\int_\Omega p  \diver\bv = \int_\Omega \bm{f}\cdot\bv\quad & \forall\,\bv\in C^\infty_0(\Omega)^d,\\
	(\BD(&\bu),\BS)\in \mathcal{A}(\cdot) & \textrm{a.e. in }\Omega,\\
	-\int_\Omega &q  \diver\bu = 0 & \forall\, q\in C^\infty_0(\Omega).
\end{alignat}
\end{subequations}
Under these conditions the results of Bul\'i\v{c}ek et al.~\cite{Bulicek:2009,Bulicek:2012} guarantee the existence of a weak solution  to problem \eqref{eq:PDE_continuous}. The results from \cite{Bulicek:2009,Bulicek:2012} apply to a class of constitutive relations more general than \eqref{eq:ImplicitCR}, but a large proportion of the algebraic constitutive relations found in practice are of the kind \eqref{eq:ImplicitCR}; in addition, when using relations of the form \eqref{eq:ImplicitCR} the conditions (C1)--(C4) are an immediate consequence of (A1)--(A4).

\begin{figure}%
	\centering
	\subfloat[A][{\centering $\nu=1$, $r_1 = 1.2$, $r_2=2.5$, \par $\beta_1=0.5$, $\beta_2=0$, $\Gamma_1=\Gamma_2=10$.}]{{
		\begin{tikzpicture}
			\begin{axis}[
				faceted color= violet,
				colormap/violet,
				grid=both,
				title={$\mu_{\textrm{eff}}(|\BS|,|\BD|)$},
			xlabel=$|\BD|$,ylabel=$|\BS|$,small,
		]
			\addplot3 [
				samples = 15,
				surf,
				domain=0:5,
				domain y=0:5,
				] {((0.5 + 0.5*pow(1 + 10*x*x,(1.2-2)/2))/(pow(1 + 10*y*y,(2-2.5)/(2*(2.5-1))))
		)};
		\end{axis}
		\end{tikzpicture}
	}}%
	\subfloat[B][{\centering $\nu=1$, $r_1 = 2.5$, $r_2=3.5$,\par $\beta_1=\beta_2=0$, $\Gamma_1=100$, $\Gamma_2=20$.}]{{
		\begin{tikzpicture}
			\begin{axis}[
				faceted color= violet,
				colormap/violet,
				grid=both,
				title={$\mu_{\textrm{eff}}(|\BS|,|\BD|)$},
			xlabel=$|\BD|$,ylabel=$|\BS|$,small,
		]
			\addplot3 [
				samples = 15,
				surf,
				domain=0:5,
				domain y=0:5,
				] {(pow(1 + 100*x*x,(2.5-2)/2))/(pow(1 + 20*y*y,(2-3.5)/(2*(3.5-1))))
		};
		\end{axis}
		\end{tikzpicture}
	}}%
	\caption{Effective viscosity for the generalised Carreau--Yasuda relation \eqref{eq:stress-Carreau}. Shear-thinning and stress-thickening behaviour can be observed in (a), while (b) presents only thickening behaviour.}%
	\label{fig:eff_visc}
\end{figure}
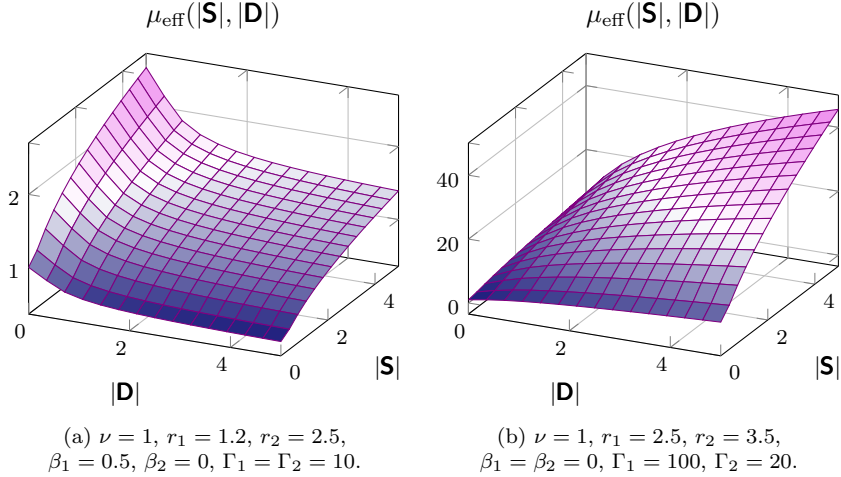

The relation \eqref{eq:ImplicitCR} defines a general constitutive law with a power law structure describing a fluid with an effective viscosity that depends both on $|\BD|$ and $|\BS|$; in this setting the effective viscosity can be defined as:
\begin{equation}
	\mu_{\text{eff}}(\cdot,|\BS|,|\BD|) := \frac{1}{2}\frac{\alpha(\cdot,|\BS|^2,|\BD|^2)}{\beta(\cdot,\abs{\BS}^2,\abs{\BD}^2)}.
\end{equation}
An important example that is captured by the assumptions above is the generalised Carreau--Yasuda \cite{Yasuda1979} constitutive relation:
\begin{equation}\label{eq:stress-Carreau}
	\begin{split}
		\BG(\BS,\BD)& =	\left(\beta_1 + (1-\beta_1)(1 + \Gamma_1\abs{\BD}^2)^{\frac{r_1-2}{2}}\right)\BD \\
	&- \frac{1}{2 \nu}\left(\beta_2 + (1-\beta_2)(1 + \Gamma_2\abs{\BS}^2)^{\frac{2-r_2}{2(r_2-1)}}\right)\BS,
	\end{split}
\end{equation}
where $r_1,r_2>1,$ $1\geq\beta_1,\beta_2\geq0$ and $\nu,\Gamma_1,\Gamma_2>0$ are given parameters.  Note that when $\beta_2 = 1$ the relation \eqref{eq:stress-Carreau} reduces to the Carreau--Yasuda constitutive relation and when $r_1=2=r_2$ or $\beta_1 = 1 = \beta_2$ it reduces to the usual Newtonian relation $\BS = 2\nu\BD$. Examples of materials that can be modelled by relations of this type when $r_2 = 2$ include molten polystyrene, ball-point pen ink, and blood, among others (see e.g.\ \cite{Barnes1989,Yasuda1979,Abraham2005}); when $r_1=2$, i.e.\ when dealing with stress-dependent viscosities, examples include ice, poly(vinyl chloride) solutions and molten polyethylene \cite{Matsuhisa1965,Savvas1994,Glen1955,Pettit2003}. \cref{fig:eff_visc} shows the behaviour of the effective viscosity for two choices of the parameters.

Another example is given by regularizations of the Bingham constitutive relation for viscoplastic fluids, which is defined by:
\begin{equation}\label{eq:Bingham}
\renewcommand{\arraystretch}{1.5}
\left\{
\begin{array}{cc}
\BS = \tau_y\frac{\BD}{|\BD|} + 2\nu \BD, & \textrm{ if } |\BS|\geq \tau_y,\\
  \BD = 0, & \textrm{ if }|\BS|<\tau_y,\\
\end{array}
\right.
\end{equation}
where $\nu>0$ and $\tau_y\geq 0$. Some examples of materials modelled by \eqref{eq:Bingham} or its power-law-like generalisation, the Herschel--Bulkley constitutive relation, include waxy crude oil, paint, pastes, drilling muds,  and mango jam \cite{Bird1983,Glowinski2010,Basu}. Note that such a relation can be written using an expression of the form \eqref{eq:ImplicitCR}; for instance, it could be described using the following functions:

\begin{subequations}\label{eq:ImplicitBingham}
	\begin{align}
		\BG_1(\BS,\BD) &= 2\nu(\tau_y + |2\nu\BD|)\BD - |2\nu\BD|\BS,\\
	\BG_2(\BS,\BD) &=\left\{
		\begin{array}{cc}
			\BD - \frac{1}{2\nu}(|\BS| - \tau_y)^+\frac{\BS}{|\BS|}, & \textrm{ if } \BS\neq \b0,\\
			\BD, & \textrm{ if } \BS = 0.\\
	\end{array}
	\right.
	\end{align}
\end{subequations}
However, the expressions in \eqref{eq:ImplicitBingham} do not satisfy the differentiability assumption (A1) and so Newton's method cannot be directly applied. This difficulty can be addressed by applying a suitable regularization step. For example, the following functions could be used instead of \eqref{eq:ImplicitBingham}:
\begin{subequations}\label{eq:ImplicitBingham_Reg}
	\begin{align}
	\tilde{\BG}_1(\BS,\BD) &= 2\nu(\tau_y + \sqrt{|2\nu\BD|^2 + \varepsilon^2})\BD - \sqrt{|2\nu\BD|^2 + \varepsilon^2}\BS,\label{eq:ImplicitBingham_Reg1}\\
	\tilde{\BG}_2(\BS,\BD) &= \left(2\nu + \frac{\tau_y}{|\bm{D}|}\right)(1- \e^{-|\BD|/\varepsilon})\BD - \BS,\label{eq:ImplicitBingham_Reg2}
\end{align}
	\end{subequations}
	where $\varepsilon$ is a positive small parameter. The relation defined by \eqref{eq:ImplicitBingham_Reg2} is known as the Papanastasiou regularization and is widely used in the simulation of viscoplastic flow \cite{Papanastasiou1987}, and while \eqref{eq:ImplicitBingham_Reg1} is related to the Bercovier--Engelman regularization \cite{Bercovier1980}, it is not usually written in this manner. This illustrates the wide freedom that the framework employed here offers; the practitioner may select the most convenient expression for a given constitutive relation.

\section{Finite element discretization}
For a barycentrically refined triangulation $\mathcal{T}_h$ of $\Omega$, let us introduce the following finite element spaces for $k\geq d$:
\begin{align*}
	\Sigma^h &= \{\bsigma\in L^\infty(\Omega)^{d\times d}\, :\, \bsigma|_K\in\mathbb{P}_{k-1}(K)^{d\times d}\text{ for all }K\in \mathcal{T}_h\},\\
V^h &= \{\bw\in W^{1,r}(\Omega)^d\, :\,\bw|_{\partial\Omega} = \bu_0,\, \bw|_K\in\mathbb{P}_{k}(K)^d\text{ for all }K\in \mathcal{T}_h\},\\
M^h &= \{q\in L_0^\infty(\Omega)\, :\, q|_K\in \mathbb{P}_{k-1}(K)\text{ for all }K\in \mathcal{T}_h\}.
\end{align*}
Here $\mathbb{P}_{k}(K)$ denotes the space of polynomials on $K$ of degree at most $k$. The velocity-pressure pair $V^h$--$M^h$ is commonly known as the Scott--Vogelius element and it is known to be inf-sup stable on barycentrically refined meshes \cite{Qin1994,Zhang2005,Tscherpel2018}, i.e.\ there is a positive constant $\alpha_1$, independent of $h$, such that
\begin{equation}
	\adjustlimits \inf_{q\in M^h\setminus\{0\}}\sup_{\bv\in V^h\setminus\{0\}} \frac{\int_\Omega q\,\diver\bv}{\|\bv\|_{W^{1,\tilde{r}'}(\Omega)}\|q\|_{L^{\tilde{r}}(\Omega)}} \geq \alpha_1.\label{infsupVel_A}\\
\end{equation}
The inf-sup condition is written in terms of $L^{\tilde{r}}(\Omega)$-norms because in this form \eqref{infsupVel_A} yields the necessary estimates for the theoretical convergence analysis; see e.g.\ \cite{Diening:2013} for details. Furthermore, we know that there exists a constant $c_r>0$ such that for every $\bsigma\in\Sigma_{\sym,\tr}^h$ there is $\btau\in \Sigma_{\sym,\tr}^h$ such that \cite[Proposition 3.1]{Sandri1998}:
\begin{equation*}
\int_\Omega \btau:\bsigma = \|\bsigma\|_{L^r(\Omega)}^r\quad\text{ and   }\quad \|\btau\|_{L^{r'}(\Omega)}\leq c_r\|\bsigma\|_{L^r(\Omega)}^{r-1}.
\end{equation*}
Combining this with the fact that $\BD(V^h)\subset \Sigma^h$, we immediately obtain inf-sup stability for the stress-velocity pair, i.e.\ there is a constant $\alpha_2>0$, independent of $h$, such that
\begin{equation}
	\adjustlimits\inf_{\bv\in V^h_{\diver}\setminus\{0\}}\sup_{\btau\in \Sigma^h_{\sym,\tr}\setminus\{0\}}  
\frac{\int_\Omega \bm{\tau}:\BD(\bm{v})}{\|\bm{\tau}\|_{L^{r'}(\Omega)}\|\bm{v}\|_{W^{1,r}(\Omega)}} \geq \alpha_2,\label{infsupStress_A}
\end{equation}
where $V^h_{\diver}$  denotes the subspace of discretely divergence-free functions of $V^h$ and $\Sigma^h_{\sym,\tr}$ denotes the subspace of symmetric and traceless functions of $\Sigma^h$. The fact that we can work with traceless stresses, and thus fewer degrees of freedom, stems from the fact that the discretely divergence-free velocities with the Scott--Vogelius element are in fact pointwise divergence-free.  This property is highly desirable and its importance has been recognized in recent years; see for instance the discussion in \cite{John2017,Linke2019}. 

In the finite element formulation of \eqref{eq:PDE_continuous} we look for $(\BS,\bu,p)\in \Sigma^h_{\sym,\tr}\times V^h\times M^h$ such that 
\begin{subequations}\label{eq:FEM_nonlinear}
\begin{alignat}{2}
	\int_\Omega \BG(\cdot &, \BS,\BD(\bu)):\btau = 0 & \forall\,\btau\in\Sigma^h_{\sym,\tr},\\
	\int_\Omega \BS:\BD(\bv) -\int_\Omega \bu\otimes\bu&:\BD(\bv) -\int_\Omega p\diver\bv = \int_\Omega \bm{f}\cdot\bv\quad & \forall\,\bv\in V^h, \label{eq:FEM_nonlinear_momentum}\\
	-&\int_\Omega q  \diver\bu = 0 & \forall\, q\in M^h.\label{eq:discrete_divergence_condition}
\end{alignat}
\end{subequations}
Known results guarantee the existence of solutions of this finite element discretization and that they each converge to a weak solution of \eqref{eq:PDE_continuous} as the mesh is refined in both the stationary \cite{Diening:2013} and transient \cite{Farrell2019} cases.

The nonlinear system \eqref{eq:FEM_nonlinear} is solved using Newton's method. Denoting the current guess for the solution as $(\tilde{\BS},\tilde{\bu},\tilde{p})$, the solution procedure is defined by a correction step $(\tilde{\BS},\tilde{\bu},\tilde{p})\mapsto(\tilde{\BS},\tilde{\bu},\tilde{p}) + (\BS,\bu,p)$ that is applied iteratively, where $(\BS,\bu,p)$ is computed by solving a linear system\textcolor{black}{, whose associated matrix presents} the following block structure:
\begin{equation}\label{eq:block_structure}
	\begin{bmatrix}
		Q_1&Q_2C^{\top}&0\\
		C&E&\tilde{B}^\top\\
		0&\tilde{B}&0
 \end{bmatrix}\begin{bmatrix}
 \BS\\ \bu\\ p
 \end{bmatrix}.
\end{equation}
The linear operators in the matrix above are defined through the relations
\begin{subequations}
	\begin{alignat}{2}
		\langle Q_1 \bsigma,\btau\rangle &:= \int_\Omega\partial_{\BS}\BG(\cdot,\tilde{\BS},\BD(\tilde{\bu})) \bsigma:\btau &\forall\,\bsigma,\btau\in \Sigma^h_{\sym,\tr},\\
		\langle Q_2C^\top \bv,\btau\rangle &:= \int_\Omega\partial_{\BD}\BG(\cdot,\tilde{\BS},\BD(\tilde{\bu})) \BD(\bv):\btau\quad &\forall\,\bv\in V^h,\btau\in \Sigma^h_{\sym,\tr},\\
		\langle & \tilde{B}  \bv,q \rangle := -\int_\Omega q\diver\bv &\forall\,\bv\in V^h,\, q\in M^h,\\
		\langle & C \bsigma,\bw  \rangle := \int_\Omega \bsigma:\BD(\bw) &\forall\,\bsigma\in \Sigma^h_{\sym,\tr},\bw\in V^h,\\
		\langle E\bv,\bw \rangle &:= -\int_\Omega (\tilde{\bu}\otimes \bv + \bv\otimes\tilde{\bu}):\BD(\bw) &\forall\,\bv,\bw\in V^h.
	\end{alignat}
\end{subequations}

Note that the differentiability and the strict monotonicity imply, together with the Implicit Function Theorem, that $Q_2^{-1}Q_1$ is either positive or negative definite. \textcolor{black}{If the convective term is neglected (or if Picard linearization is used instead), with the help of the inf-sup conditions \eqref{infsupVel_A} and \eqref{infsupStress_A} one can guarantee that \eqref{eq:block_structure} is invertible. Although the invertibility of \eqref{eq:block_structure} is not clear when using Newton's method, in this work we will always employ it, because of its quadratic convergence rate (assuming the current guess is sufficiently close to the solution).}

\section{Augmented Lagrangian Preconditioner}\label{Sec:AL}
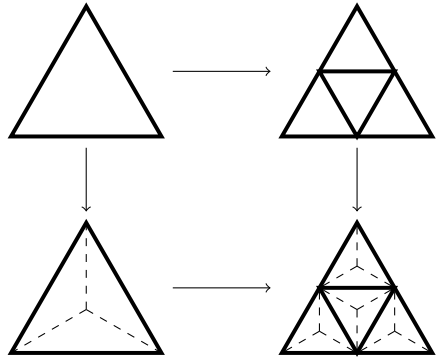
\begin{figure}
\begin{center}
\begin{tikzpicture}[scale=1.8]
\usetikzlibrary{positioning}
	\node (A) [] at (-1,1){
		\begin{tikzpicture}[scale=2]
			\draw[line width=1.5pt] (-0.5,-0) -- (0.5,0) -- (0.,0.866) -- cycle;
		\end{tikzpicture}
		};
	\node (B) [] at (1,1) {
		\begin{tikzpicture}[scale=2]
			\draw[line width=1.5pt] (-0.5,-0) -- (0.5,0) -- (0.,0.866) -- cycle;
			\draw[line width=1.5pt] (-0,-0) -- (0.25,0.433) -- (-0.25,0.433) -- cycle;
		\end{tikzpicture}
		};
	\draw [->] (A) -- (B);
	\node (A_bary) [] at (-1,-0.6){
		\begin{tikzpicture}[scale=2]
			\draw[line width=1.5pt] (-0.5,-0) -- (0.5,0) -- (0.,0.866) -- cycle;
			\begin{scope}[style=dashed]
				\coordinate (bar1) at (-0.5,0);
				\coordinate (bar2) at (0.5,0);
				\coordinate (bar3) at (0,0.866);
				\draw (barycentric cs:bar1=0.5,bar2=0.5,bar3=0.5) -- (barycentric cs:bar1=1,bar2=0,bar3=0);
				\draw (barycentric cs:bar1=0.5,bar2=0.5,bar3=0.5) -- (barycentric cs:bar1=0,bar2=1,bar3=0);
				\draw (barycentric cs:bar1=0.5,bar2=0.5,bar3=0.5) -- (barycentric cs:bar1=0,bar2=0,bar3=1);
			\end{scope}
		\end{tikzpicture}
		};
	\node (B_bary) [] at (1,-0.6) {
		\begin{tikzpicture}[scale=2]
			\draw[line width=1.5pt] (-0.5,-0) -- (0.5,0) -- (0.,0.866) -- cycle;
			\draw[line width=1.5pt] (-0,-0) -- (0.25,0.433) -- (-0.25,0.433) -- cycle;
			\begin{scope}[style=dashed]
			\foreach \x in {-0.5,0.}{
				\coordinate (bar1) at (\x,0);
				\coordinate (bar2) at (\x+0.5,0);
				\coordinate (bar3) at (\x+0.25,0.433);
				\draw (barycentric cs:bar1=0.5,bar2=0.5,bar3=0.5) -- (barycentric cs:bar1=1,bar2=0,bar3=0);
				\draw (barycentric cs:bar1=0.5,bar2=0.5,bar3=0.5) -- (barycentric cs:bar1=0,bar2=1,bar3=0);
				\draw (barycentric cs:bar1=0.5,bar2=0.5,bar3=0.5) -- (barycentric cs:bar1=0,bar2=0,bar3=1);
			}
				\coordinate (bar1) at (0.,0);
				\coordinate (bar2) at (-0.25,0.433);
				\coordinate (bar3) at (0.25,0.433);
				\draw (barycentric cs:bar1=0.5,bar2=0.5,bar3=0.5) -- (barycentric cs:bar1=1,bar2=0,bar3=0);
				\draw (barycentric cs:bar1=0.5,bar2=0.5,bar3=0.5) -- (barycentric cs:bar1=0,bar2=1,bar3=0);
				\draw (barycentric cs:bar1=0.5,bar2=0.5,bar3=0.5) -- (barycentric cs:bar1=0,bar2=0,bar3=1);
				\coordinate (bar1) at (0,0.866);
				\coordinate (bar2) at (-0.25,0.433);
				\coordinate (bar3) at (0.25,0.433);
				\draw (barycentric cs:bar1=0.5,bar2=0.5,bar3=0.5) -- (barycentric cs:bar1=1,bar2=0,bar3=0);
				\draw (barycentric cs:bar1=0.5,bar2=0.5,bar3=0.5) -- (barycentric cs:bar1=0,bar2=1,bar3=0);
				\draw (barycentric cs:bar1=0.5,bar2=0.5,bar3=0.5) -- (barycentric cs:bar1=0,bar2=0,bar3=1);
			\end{scope}
		\end{tikzpicture}
		};
	\draw [->] (A_bary) -- (B_bary);
	\draw [->] (A) -- (A_bary);
	\draw [->] (B) -- (B_bary);
\end{tikzpicture}
\end{center}
\caption{\label{fig:mesh_hierarchy}
Non-nested two-level barycentrically refined mesh hierarchy \textcolor{black}{in two dimensions}.}
\end{figure}

As mentioned above, barycentric refinement guarantees the inf-sup stability of the Scott--Vogelius element pair for $k \ge d$.
However, constructing a multigrid hierarchy by successive barycentric refinement creates degenerate elements.
We therefore employ the alternative construction used in \cite{Farrell2020}. The multigrid hierarchy is obtained by taking a standard uniformly-refined hierarchy and barycentrically refining on each level once; see
\cref{fig:mesh_hierarchy}. The cells before barycentric refinement are referred to as \emph{macro cells}.
An important consequence of this is the existence of local Fortin operators on each macro cell, which are useful when trying to characterize locally the space of divergence-free velocities \cite{Farrella}. 
A disadvantage is that the resulting mesh hierarchy is non-nested, which leads to some complications with the prolongation operator in the multigrid algorithm.

\begin{figure}
\begin{center}
\begin{tikzpicture}[scale=1.5]
\usetikzlibrary{shapes.geometric}
\begin{scope}
\draw[step=1cm] (-2,-1) grid (2,2);
\draw (-2,2) -- (1,-1);
\draw (-1,2) -- (2,-1);
\draw (0,2) -- (2,0);
\draw (1,2) -- (2,1);
\draw (-2,0) -- (-1,-1);
\draw (-2,-1) -- (-1,-1);
\draw (-2,1) -- (0,-1);
\end{scope}

\begin{scope}[style=help lines]
\foreach \x  in {0,1,-1,-2} {
\foreach \y  in {0,1,-1,-1} {
\coordinate (bar)   at (\x,\y);
\coordinate (barr) at (\x,\y + 1);
\coordinate (barrr)      at (\x + 1,\y);
\coordinate (cbar)   at (\x + 1,\y + 1);
\draw (barycentric cs:cbar=0.5,barr=0.5 ,barrr=0.5) -- (barycentric cs:cbar=1.,barr=0. ,barrr=0);
\draw (barycentric cs:cbar=0.5,barr=0.5 ,barrr=0.5) -- (barycentric cs:cbar=0.,barr=1 ,barrr=0);
\draw (barycentric cs:cbar=0.5,barr=0.5 ,barrr=0.5) -- (barycentric cs:cbar=0.,barr=0. ,barrr=1);
\draw (barycentric cs:bar=0.5,barr=0.5 ,barrr=0.5) -- (barycentric cs:bar=1.,barr=0. ,barrr=0);
\draw (barycentric cs:bar=0.5,barr=0.5 ,barrr=0.5) -- (barycentric cs:bar=0.,barr=1 ,barrr=0);
\draw (barycentric cs:bar=0.5,barr=0.5 ,barrr=0.5) -- (barycentric cs:bar=0.,barr=0. ,barrr=1);}
}
\end{scope}

\node[regular polygon, regular polygon sides=6, draw] [dashed, ultra thick, rounded corners=8pt, rotate=-45,xscale=10.5,yscale=6,magenta] at (-1,0){}; 
\filldraw [magenta] (-1,0) circle (1.5pt);
\node[regular polygon, regular polygon sides=6, draw] [dashed, ultra thick,rounded corners=8pt,rotate=-45,xscale=10.5,yscale=6.,magenta] at (1,1){};
\filldraw [magenta] (1,1) circle (1.5pt);
\end{tikzpicture}
\end{center}
\caption{\label{fig:macrostar}
Macrostar patches on a barycentrically refined mesh \textcolor{black}{in two dimensions}.}
\end{figure}
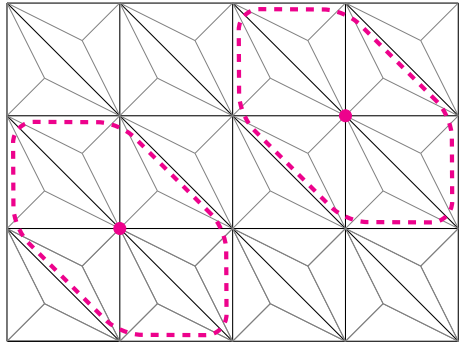

\begin{remark}
Augmented Lagrangian preconditioners have been applied to flow problems with variable viscosity and Bingham rheology before; see e.g.\ \cite{He2015,He2012}. In those works it is advocated that for the Schur complement approximation a viscosity-weighted mass matrix should be used instead:
\begin{equation}\label{eq:variable_mass_matrix}
	(M_\mu)_{ij} := \int_\Omega \frac{1}{\mu}\phi_i\,\phi_j,
\end{equation}
where $\mu$ denotes the variable (effective) viscosity and $\phi_i, \phi_j$ are pressure basis functions.
A similar argument was presented in \cite{Grinevich2009}, where only the Schur complement approximation without the augmented Lagrangian term was studied. However, in those works a robust scalable solver for the augmented momentum block was not available and the authors were limited to low values of $\gamma$ ($\gamma=1$ was used in their numerical experiments), and so a better approximation for the Schur complement with \eqref{eq:variable_mass_matrix} was necessary. In contrast, the multigrid solver presented in this work for the \textcolor{black}{stress-velocity} block will be $\gamma$-robust, which therefore allows for very large values of $\gamma$, and thus excellent control of the Schur complement. It is consequently not necessary to use \eqref{eq:variable_mass_matrix}, which requires reassembly at every Newton step.

\subsubsection*{Robust Relaxation}
\textcolor{black}{From \eqref{eq:AL_stabilisation} and \eqref{eq:block_structure}, we see that the augmented stress-velocity block can be written as}
\begin{equation}\label{eq:top_block}
	A + \gamma B^{\top}M_p^{-1}B =
	\begin{bmatrix}
		Q_1 & Q_2 C^\top \\ C & E
	\end{bmatrix}
	+ \gamma \begin{bmatrix}
		0 \\ \tilde{B}^\top 
	\end{bmatrix}
	M_p^{-1}\begin{bmatrix}
		0 & \tilde{B}
	\end{bmatrix},
\end{equation}
where $\gamma B^{\top}M_p^{-1}B$ is symmetric and semidefinite and $A$ is invertible. Relaxation methods used in multigrid algorithms can often be framed in terms of subspace correction methods \cite{Xu1992,Xu2001}. Let us define $Z^h := \Sigma_{\sym,\tr}^h\times V^h$ and consider the space decomposition
\begin{equation}\label{eq:subspace_decomposition}
Z^h = \sum_i Z_i,
\end{equation}
where the sum is not necessarily direct; for instance, Jacobi relaxation is obtained by setting $Z_i = \lin \{\bphi_i\}$, where $\{\bphi_i\}_i$ denotes a basis of $Z$. For the case when $A$ is symmetric and positive definite the theory is rigorously established; the results of Sch\"oberl~\cite{Schoeberl1999} and Lee et al.~\cite{LWXZ:2007} guarantee the $\gamma$-robustness of the subspace correction method provided that the decomposition stably captures the kernel of the semidefinite term. Capturing the kernel means that
\begin{equation}
	\mathcal{N}^h = \sum_i Z_i\cap \mathcal{N}^h,
\end{equation}
where $\mathcal{N}^h$ denotes the kernel of the semidefinite term (in the case of \eqref{eq:top_block}, these are the elements of the form $(\btau, \bv)^\top$, where $\bv$ is divergence-free in $V^h$, and $\btau\in \Sigma^h_{\sym,\tr}$ is arbitrary).

Thankfully, a local characterisation of the kernel of the divergence operator for the Scott--Vogelius discretization on meshes with the \textcolor{black}{macro structure considered here}, was recently obtained in \cite{Farrella} (see also \cite{Farrell2020}). In that work it was proven that a kernel capturing space decomposition is obtained by setting 
\begin{equation}\label{eq:decomposition_SV}
	Z_i := \{\bz\in Z^h\, :\, \supp(\bz)\subset \macrostar(q_i)  \},
\end{equation}
where for each vertex $q_i$, the macrostar patch $\macrostar(q_i)$ is defined as the union of all macro cells touching the vertex (\cref{fig:macrostar} \textcolor{black}{shows a two-dimensional example}).
\end{remark}

\begin{remark}\label{rem:ActEuler}
	In some cases the analysis can be carried out in a slightly different manner. For example, if we take a Bercovier--Engelman-like regularization of the constitutive relation for an activated Euler fluid (this is the counterpart of the Bingham constitutive relation where the roles of $\BS$ and $\BD$ are interchanged, see e.g.\ \cite{Blechta2019})
	\begin{equation*}
	\BG(\BS,\BD) = \BD - \left(\frac{1}{2\nu} + \frac{\tau_y}{\sqrt{\varepsilon^2 + |\BS|^2}} \right)\BS,
\end{equation*}
with $\nu,\varepsilon>0$ and $\tau_y\geq 0$, then the \textcolor{black}{stress-velocity} block in the linearized problem can be split as follows:
\begin{equation}\label{eq:split_act_euler}
	\hat{A}_\nu + \hat{A}_\varepsilon + \gamma B^{\top}M_p^{-1}B,
\end{equation}
where $\hat{A}_\nu$ corresponds to the operator arising from the Newtonian problem and $\hat{A}_\varepsilon$ is defined via
\begin{equation*}
	\langle \hat{A}_\varepsilon (\bsigma,\bv),(\btau,\bw)\rangle := \tau_y\int_\Omega \frac{1}{\sqrt{\varepsilon^2 + |\tilde{\BS}|^2}}\left[\BI - \frac{\tilde{\BS}\otimes\tilde{\BS}}{\varepsilon^2 + |\tilde{\BS}|^2} \right]\bsigma:\btau,\quad\forall\, (\bsigma,\bv),(\btau,\bw)\in Z^h.
\end{equation*}
The splitting \eqref{eq:split_act_euler} could then be interpreted as a perturbation of the Newtonian problem, which results in an operator that degenerates as $\varepsilon\rightarrow 0$, $\gamma\rightarrow\infty$, with a kernel given by elements of the form $(\tilde{\BS},\bw)\in Z^h$, with $\diver\bw=0$. Note that while the kernel possesses a one-dimensional stress component, in practice this does not appear to cause any difficulties for the preconditioner. An illustrative example for a slightly more complicated problem will be shown in the final section of this work.
\end{remark}
In the algorithm presented here, the relaxation solves will be performed additively. For the patches depicted in \cref{fig:macrostar}, each coupled stress-velocity solve for $k=2$ (resp.~$k=3$) involves 31 (resp.~73) degrees of freedom for each component of the velocity and 60 (resp.~156) degrees of freedom for the stress. This is much more expensive than, say, a Jacobi smoother, but the resulting robustness in the algorithm makes it worth the cost, and small local patchwise solves are quite well suited to modern computing architectures.

\begin{remark}	
	When working with the full nonlinear problem including advection, the macrostar iteration \eqref{eq:subspace_decomposition} \& \eqref{eq:decomposition_SV} is not effective as a standalone relaxation method. However, as observed in \cite{Farrell2019a,Farrell2020}, this difficulty can be overcome by applying a small number of GMRES iterations preconditioned by the macrostar iteration as relaxation.
\end{remark}

\subsubsection*{Robust Prolongation}
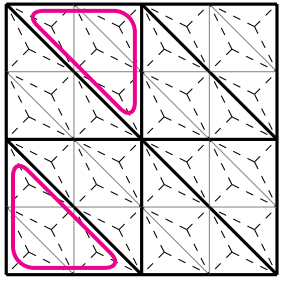
\begin{figure}
\begin{center}
\begin{tikzpicture}[scale=1.8]
\usetikzlibrary{shapes.geometric}
\begin{scope}[very thick]
\draw[step=1cm] (-1,-1) grid (1,1);
\draw (-1,1) -- (1,-1);
\draw (-1,0) -- (0,-1);
\draw (0,1) -- (1,0);
\end{scope}
\begin{scope}[style = help lines]
\foreach \x in {-1,0}{
\foreach \y in {-1,0}{
\draw (\x,\y+0.5) -- (\x+1,\y+0.5);
\draw (\x+0.5,\y) -- (\x+0.5,\y+1);
\draw (\x,\y+0.5) -- (\x+0.5,\y);
\draw (\x+0.5,\y+1) -- (\x+1,\y+0.5);
}
}
\end{scope}

\begin{scope}[style=dashed]
\foreach \x  in {0.5,0,-0.5,-1} {
\foreach \y  in {0.5,0,-0.5,-1} {
\coordinate (bar)   at (\x,\y);
\coordinate (abar) at (\x+0.5,\y+0.5);
\coordinate (bbar)      at (\x,\y+0.5);
\coordinate (cbar)   at (\x+0.5,\y);
\draw (barycentric cs:bar=0.5,bbar=0.5 ,cbar=0.5) -- (barycentric cs:bar=1.,bbar=0. ,cbar=0);
\draw (barycentric cs:bar=0.5,bbar=0.5 ,cbar=0.5) -- (barycentric cs:bar=0.,bbar=1. ,cbar=0);
\draw (barycentric cs:bar=0.5,bbar=0.5 ,cbar=0.5) -- (barycentric cs:bar=0.,bbar=0. ,cbar=1);
\draw (barycentric cs:abar=0.5,bbar=0.5 ,cbar=0.5) -- (barycentric cs:abar=1.,bbar=0. ,cbar=0);
\draw (barycentric cs:abar=0.5,bbar=0.5 ,cbar=0.5) -- (barycentric cs:abar=0.,bbar=1. ,cbar=0);
\draw (barycentric cs:abar=0.5,bbar=0.5 ,cbar=0.5) -- (barycentric cs:abar=0.,bbar=0. ,cbar=1);
}
}
\end{scope}
%
\draw[line width=1.5pt,color=magenta,rounded corners=8pt] (-0.95,-0.95) -- (-0.12,-0.95) -- (-0.95,-0.12) -- cycle;
\draw[line width=1.5pt,color=magenta,rounded corners=8pt] (-0.88,0.95) -- (-0.05,0.95) -- (-0.05,0.12) -- cycle;
\end{tikzpicture}
\end{center}
	\caption{\label{fig:patch_prolongation}
	The correction to the prolongation operator is computed on the coarse macro cells.}
\end{figure}

A robust multigrid algorithm also requires a stable prolongation operator $P_H:Z^H\rightarrow Z^h$, mapping the space of coarse grid functions $Z^H$ into the space of fine grid functions $Z^h$, with a continuity constant independent of $\gamma$ or parameters arising in the implicit constitutive relation. In the setting of a velocity-pressure formulation of the Stokes problem, the matrix $A$ acts only on the velocity space $V^h$ and is actually SPD and thus the whole matrix \eqref{eq:top_block} defines a norm. We could therefore write:
\begin{align*}
	\|\bv_H\|_{H,\gamma}^2 &= \|\bv_H\|^2_{A_H} + \gamma \|\diver\bv_H\|_{L^2(\Omega)}^2,\\
	\|P_H\bv_H\|_{h,\gamma}^2 &= \|P_H\bv_H\|^2_{A_h} + \gamma \|\diver(P_H\bv_H)\|_{L^2(\Omega)}^2,
\end{align*}
where $A_H$ and $A_h$ correspond to discretizations on the coarse and fine mesh, respectively. The central difficulty is that the condition $\diver\bv_H=0$ does not necessarily imply that $\diver(P_H\bv_H) = 0$, when $P_H$ is a standard prolongation operator based on finite element interpolation, due to the non-nestedness of the mesh hierarchy. If not addressed, this causes a lack of robustness in the multigrid solver for large $\gamma$. The insight of Sch\"{o}berl \cite{Schoeberl:1999,Schoeberl1999}, later applied by Benzi and Olshanskii in \cite{Benzi2006}, and Farrell, Mitchell and Wechsung \cite{Farrell2019a,Farrella,Farrell2020}, is that by performing local Stokes solves it is possible to compute a correction to the prolongation operator and ensure that divergence-free fields get mapped to (nearly) divergence-free fields. For the Scott--Vogelius discretization on meshes with the macro structure illustrated in  \cref{fig:mesh_hierarchy}, it can be seen that interpolation is actually exact on the boundaries of the coarse macro cells, and therefore, as shown in \cite{Farrella,Farrell2020}, the correction to the prolongation operator can be computed on the space
\begin{equation*}
	\tilde{V}^h := \{\bv_h \in V^h \, :\, \supp(\bv)\subset K\text{ for some }K\in\mathcal{M}_H\},
\end{equation*}
where $\mathcal{M}_H$ is the triangulation of coarse macro elements. To be more precise, the corrected prolongation operator is defined through
\begin{equation}
	\tilde{P}_H \bv_H := P_H\bv_H - \tilde{\bv}_h,
\end{equation}
where $\tilde{\bv}_h \in \tilde{V}^h$ solves the Stokes-like problem
\begin{equation}\label{eq:local_solve_prolongation}
	\int_\Omega 2\nu \BD(\tilde{\bv}_h):\BD(\tilde{\bw}_h) + \gamma\int_\Omega \diver\tilde{\bv}_h\,\diver\tilde{\bw}_h = \gamma \int_\Omega \diver(P_H\bv_H)\,\diver\tilde{\bw}_h\quad \forall\,\tilde{\bw}_h\in \tilde{V}^h.
\end{equation}
The positive parameter $\nu$ is arbitrary; it could for instance be taken as the one appearing in the constitutive relations \eqref{eq:stress-Carreau} and \eqref{eq:Bingham}. Observe that, by definition of the space $\tilde{V}^h$, the problem \eqref{eq:local_solve_prolongation} decouples on the patches defined by the macro elements and can therefore be computed independently on each macro cell (see  \cref{fig:patch_prolongation}); this is important for the efficiency of the solver.

In the non-Newtonian setting, it may seem more appropriate to alternatively employ on the left hand side of \eqref{eq:local_solve_prolongation} the operator defined by the Schur complement $-CQ_1^{-1}Q_2C^{\top}$ (which reduces to \eqref{eq:local_solve_prolongation} in the Newtonian case). However, since the end goal is to correct for the error in the divergence introduced by the interpolation operator, we prefer to retain \eqref{eq:local_solve_prolongation} for the sake of avoiding reassembly and refactorization.

The prolongation operator for the stress variables $\bsigma_H \in \Sigma^H \mapsto \bsigma_h\in \Sigma^h$, between spaces $\Sigma^H$ and $\Sigma^h$ defined on the coarse and fine meshes respectively, is defined via the Galerkin projection
		\begin{equation}\label{eq:Galerkin_projection}
			\|\bsigma_h - \bsigma_H\|_{L^2(\Omega)} =  \min_{\bsigma\in \Sigma^h}\|\bsigma_H - \bsigma\|_{L^2(\Omega)}.
		\end{equation}
If we denote the basis of $\Sigma^h$ by $\{\bphi_h^i\}_{i=1}^{N_h}$, then the optimality condition for \eqref{eq:Galerkin_projection} takes the form
		\begin{equation} \label{eq:GPoptimality}
            \int_\Omega \bsigma_h\colon \bphi_h^i = \int_\Omega \bsigma_H \colon \bphi_h^i\quad \forall\, i\in\{1,\ldots,N_h\},
		\end{equation}
or written in matrix form:
\begin{equation}
	M_h\bsigma_h = M_{h,H}\bsigma_H,
\end{equation}
where the mass matrices are defined as
\begin{equation}\label{eq:mixed_mass_matrices}
	\begin{gathered}
		(M_h)_{ij} = \int_\Omega \bphi^i_h\colon \bphi^j_h\quad i,j\in\{1,\ldots,N_h\},\\
		(M_{h,H})_{ij} = \int_\Omega \bphi^i_h\colon \bphi^j_H\quad i\in\{1,\ldots,N_h\},\, j\in\{1,\ldots,N_H\},\\
	\end{gathered}
\end{equation}
where the basis of $\Sigma^H$ is denoted by $\{\bphi_H^i\}_{i=1}^{N_H}$.

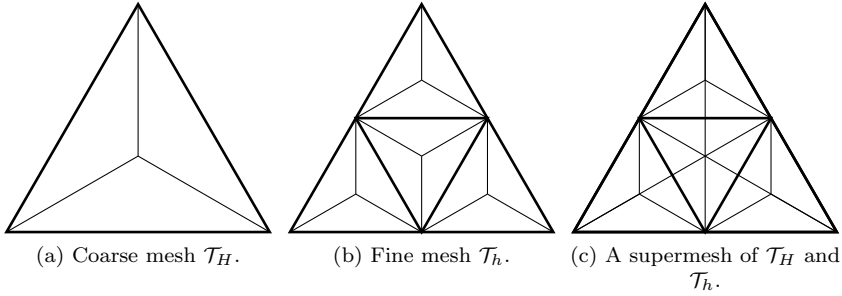
\begin{figure}%
	\centering
	\subfloat[A][{\centering Coarse mesh $\mathcal{T}_H$.}]{{
		\begin{tikzpicture}[scale=3.5]
			\draw[line width=1pt] (-0.5,-0) -- (0.5,0) -- (0.,0.866) -- cycle;
			\begin{scope}[]
				\coordinate (bar1) at (-0.5,0);
				\coordinate (bar2) at (0.5,0);
				\coordinate (bar3) at (0,0.866);
				\draw (barycentric cs:bar1=0.5,bar2=0.5,bar3=0.5) -- (barycentric cs:bar1=1,bar2=0,bar3=0);
				\draw (barycentric cs:bar1=0.5,bar2=0.5,bar3=0.5) -- (barycentric cs:bar1=0,bar2=1,bar3=0);
				\draw (barycentric cs:bar1=0.5,bar2=0.5,bar3=0.5) -- (barycentric cs:bar1=0,bar2=0,bar3=1);
			\end{scope}
		\end{tikzpicture}
	}}%
	\subfloat[B][{\centering Fine mesh $\mathcal{T}_h$.}]{{
		\begin{tikzpicture}[scale=3.5]
			\draw[line width=1pt] (-0.5,-0) -- (0.5,0) -- (0.,0.866) -- cycle;
			\draw[line width=1pt] (-0,-0) -- (0.25,0.433) -- (-0.25,0.433) -- cycle;
			\begin{scope}[]
			\foreach \x in {-0.5,0.}{
				\coordinate (bar1) at (\x,0);
				\coordinate (bar2) at (\x+0.5,0);
				\coordinate (bar3) at (\x+0.25,0.433);
				\draw (barycentric cs:bar1=0.5,bar2=0.5,bar3=0.5) -- (barycentric cs:bar1=1,bar2=0,bar3=0);
				\draw (barycentric cs:bar1=0.5,bar2=0.5,bar3=0.5) -- (barycentric cs:bar1=0,bar2=1,bar3=0);
				\draw (barycentric cs:bar1=0.5,bar2=0.5,bar3=0.5) -- (barycentric cs:bar1=0,bar2=0,bar3=1);
			}
				\coordinate (bar1) at (0.,0);
				\coordinate (bar2) at (-0.25,0.433);
				\coordinate (bar3) at (0.25,0.433);
				\draw (barycentric cs:bar1=0.5,bar2=0.5,bar3=0.5) -- (barycentric cs:bar1=1,bar2=0,bar3=0);
				\draw (barycentric cs:bar1=0.5,bar2=0.5,bar3=0.5) -- (barycentric cs:bar1=0,bar2=1,bar3=0);
				\draw (barycentric cs:bar1=0.5,bar2=0.5,bar3=0.5) -- (barycentric cs:bar1=0,bar2=0,bar3=1);
				\coordinate (bar1) at (0,0.866);
				\coordinate (bar2) at (-0.25,0.433);
				\coordinate (bar3) at (0.25,0.433);
				\draw (barycentric cs:bar1=0.5,bar2=0.5,bar3=0.5) -- (barycentric cs:bar1=1,bar2=0,bar3=0);
				\draw (barycentric cs:bar1=0.5,bar2=0.5,bar3=0.5) -- (barycentric cs:bar1=0,bar2=1,bar3=0);
				\draw (barycentric cs:bar1=0.5,bar2=0.5,bar3=0.5) -- (barycentric cs:bar1=0,bar2=0,bar3=1);
			\end{scope}
		\end{tikzpicture}
	}}%
	\subfloat[C][{\centering A supermesh of $\mathcal{T}_H$ and $\mathcal{T}_h$.}]{{
		\begin{tikzpicture}[scale=3.5]
			\draw[line width=1pt] (-0.5,-0) -- (0.5,0) -- (0.,0.866) -- cycle;
			\begin{scope}[]
				\coordinate (bar111) at (-0.5,0);
				\coordinate (bar222) at (0.5,0);
				\coordinate (bar333) at (0,0.866);
				\draw (barycentric cs:bar111=0.5,bar222=0.5,bar333=0.5) -- (barycentric cs:bar111=1,bar222=0,bar333=0);
				\draw (barycentric cs:bar111=0.5,bar222=0.5,bar333=0.5) -- (barycentric cs:bar111=0,bar222=1,bar333=0);
				\draw (barycentric cs:bar111=0.5,bar222=0.5,bar333=0.5) -- (barycentric cs:bar111=0,bar222=0,bar333=1);
			\end{scope}
			\begin{scope}[]
			\draw[line width=1pt] (-0.5,-0) -- (0.5,0) -- (0.,0.866) -- cycle;
			\draw[line width=1pt] (-0,-0) -- (0.25,0.433) -- (-0.25,0.433) -- cycle;
			\foreach \x in {-0.5,0.}{
				\coordinate (bar1) at (\x,0);
				\coordinate (bar2) at (\x+0.5,0);
				\coordinate (bar3) at (\x+0.25,0.433);
				\draw (barycentric cs:bar1=0.5,bar2=0.5,bar3=0.5) -- (barycentric cs:bar1=1,bar2=0,bar3=0);
				\draw (barycentric cs:bar1=0.5,bar2=0.5,bar3=0.5) -- (barycentric cs:bar1=0,bar2=1,bar3=0);
				\draw (barycentric cs:bar1=0.5,bar2=0.5,bar3=0.5) -- (barycentric cs:bar1=0,bar2=0,bar3=1);
			}
				\coordinate (bar1) at (0.,0);
				\coordinate (bar2) at (-0.25,0.433);
				\coordinate (bar3) at (0.25,0.433);
				\draw (barycentric cs:bar1=0.5,bar2=0.5,bar3=0.5) -- (barycentric cs:bar1=1,bar2=0,bar3=0);
				\draw (barycentric cs:bar1=0.5,bar2=0.5,bar3=0.5) -- (barycentric cs:bar1=0,bar2=1,bar3=0);
				\draw (barycentric cs:bar1=0.5,bar2=0.5,bar3=0.5) -- (barycentric cs:bar1=0,bar2=0,bar3=1);
				\coordinate (bar1) at (0,0.866);
				\coordinate (bar2) at (-0.25,0.433);
				\coordinate (bar3) at (0.25,0.433);
				\draw (barycentric cs:bar1=0.5,bar2=0.5,bar3=0.5) -- (barycentric cs:bar1=1,bar2=0,bar3=0);
				\draw (barycentric cs:bar1=0.5,bar2=0.5,bar3=0.5) -- (barycentric cs:bar1=0,bar2=1,bar3=0);
				\draw (barycentric cs:bar1=0.5,bar2=0.5,bar3=0.5) -- (barycentric cs:bar1=0,bar2=0,bar3=1);
			\end{scope}
		\end{tikzpicture}
	}}%
	\caption{Example of a coarse mesh $\mathcal{T}_H$, a fine mesh $\mathcal{T}_h$, and an associated supermesh.}%
	\label{fig:supermesh}
\end{figure}

Since the meshes are non-nested, the assembly of $M_{h,H}$ requires the integration of \emph{piecewise} polynomial functions over the cells of either mesh. To integrate these accurately we construct a \emph{supermesh} of both input meshes \cite{Farrell2011}, a common refinement of both (see \cref{fig:supermesh}). Over each supermesh cell the integrand of the right-hand side of \eqref{eq:GPoptimality} is polynomial, and hence can be calculated accurately with standard quadrature rules. Since the stress is approximated using discontinuous piecewise polynomials, the mass matrix $M_h$ is block diagonal, and is simple to invert exactly. 

As rediscretization is employed to assemble coarse grid problems, the current guess for the stress must be injected onto coarse grids. Injection is defined via a Galerkin projection analogous to \eqref{eq:Galerkin_projection}, and employs the same supermesh.

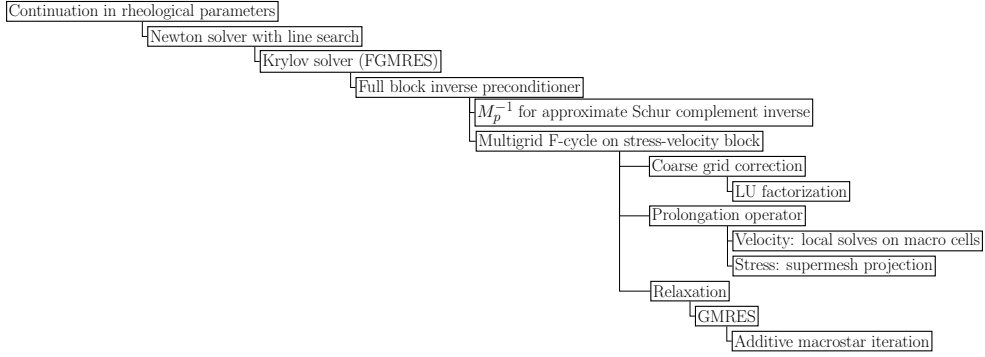
\begin{figure}
\begin{center}
\resizebox{\textwidth}{!}{%
\begin{tikzpicture}[node distance=2cm]

\usetikzlibrary{shapes.geometric, arrows}
\usetikzlibrary{positioning}

\tikzstyle{startstop} = [rectangle, minimum width=3cm, minimum height=1cm,text centered, draw=black, ultra thick]
\tikzstyle{arrow} = [thick,>=stealth]

\node (continuation) [startstop] {\Huge Continuation in rheological parameters};
\node (newton) [startstop, below=0.75cm, right=0.25cm, at=(continuation.south)] {\Huge Newton solver with line search};
\node (krylov) [startstop, below=0.75cm, right= 0.25cm, at=(newton.south)] {\Huge Krylov solver (FGMRES)};
\node (blockp) [startstop, below=0.75cm, right= 0.25cm, at=(krylov.south)] {\Huge Full block inverse preconditioner};
\node (schur) [startstop, below=0.75cm, right= 0.25cm, at=(blockp.south)] {\Huge $M^{-1}_p$ for approximate Schur complement inverse};
\node (fcycle) [startstop, below=0.75cm, right=-8.6cm, at=(schur.south)] {\Huge Multigrid F-cycle on stress-velocity block};
\node (coarseg) [startstop, below=0.75cm, right= 1.5cm, at=(fcycle.south)] {\Huge Coarse grid correction};
\node (LU) [startstop, below=0.75cm, right= 0.25cm, at=(coarseg.south)] {\Huge LU factorization};
\node (prolong) [startstop, below=2.0cm, right=-4cm, at=(coarseg.south)] {\Huge Prolongation operator};
\node (velocity) [startstop, below=0.75cm, right=0.25cm, at=(prolong.south)] {\Huge Velocity: local solves on macro cells};
\node (stress) [startstop, below=0.75cm, right=-6.4cm, at=(velocity.south)] {\Huge Stress: supermesh projection};
\node (relaxation) [startstop, below=0.75cm, right=-9.4cm, at=(stress.south)] {\Huge Relaxation};
\node (gmres) [startstop, below=0.75cm, right= 0.25cm, at=(relaxation.south)] {\Huge GMRES};
\node (macrostar) [startstop, below=0.75cm, right=0.25cm, at=(gmres.south)] {\Huge Additive {macrostar} iteration};

\draw [arrow] (continuation) |- (newton);
\draw [arrow] (newton) |- (krylov);
\draw [arrow] (krylov) |- (blockp);
\draw [arrow] (blockp) |- (schur);
\draw [arrow] (blockp) |- (fcycle);
\draw [arrow] (fcycle) |- (coarseg);
\draw [arrow] (fcycle) |- (prolong);
\draw [arrow] (fcycle) |- (relaxation);
\draw [arrow] (coarseg) |- (LU);
\draw [arrow] (prolong) |- (velocity);
\draw [arrow] (prolong) |- (stress);
\draw [arrow] (relaxation) |- (gmres);
\draw [arrow] (gmres) |- (macrostar);
\end{tikzpicture}
}
\end{center}
\caption{\label{fig:algorithm}
	Overview of the algorithm.}
\end{figure}

An overview of the full algorithm can be found in  \cref{fig:algorithm}.

\section{Numerical Examples}\label{Sec:Examples}
All the numerical examples presented in this work were \textcolor{black}{implemented} using {Firedrake} \cite{Rathgeber2016}. The macrostar patch solves for the relaxation and the local solves for the prolongation operator in the multigrid algorithm were carried out with \textsc{PCPATCH} \cite{Farrell}, a recently developed preconditioner in PETSc \cite{PETSc} for matrix-free multigrid relaxation via space decompositions. The $L^2$ line search algorithm \cite{PETScLi} was employed to improve the convergence of the Newton solver; the Newton solver was deemed to have converged when the Euclidean norm of the residual fell below $1\times 10^{-8}$ and the corresponding tolerance for the linear solver was set to $1\times 10^{-10}$\textcolor{black}{, unless specified otherwise}. \textcolor{black}{These tight tolerances are taken to challenge the solver; in practical computations the tolerance on the linear solver should be dynamically adjusted to minimise the computational work, using e.g.~the Eisenstat--Walker algorithm \cite{eisenstat1996}}. The augmented Lagrangian parameter was taken as $\gamma = 10^4$, to obtain excellent control of the Schur complement. In the implementation, the uniqueness of the pressure was recovered not by enforcing a zero mean condition in the variational formulation but rather by orthogonalizing against the nullspace of constants in the Krylov solver.
\subsubsection*{Bingham flow between two plates}
We first test our solver on a problem where the exact solution is known. Let $\Omega = (0,L)\times(-1,1)$ with $L>0$ and consider problem \eqref{eq:PDE_continuous} with $\bm{f}=\bm{0}$ and the Bingham constitutive relation \eqref{eq:Bingham}. A function that solves this problem exactly is given by \cite{Aposporidis2011,Grinevich2009,Hron2017}:
\begin{gather}\label{eq:Bingham_Poiseuille}
\bu_e(\bm{x}):= \left\{
\begin{array}{cc}
\renewcommand{\arraystretch}{1.5}
(\frac{C}{2}(1-x_2^2) - \tau_y(1-x_2),0)^\top, & \textrm{ if } \frac{\tau_y}{C}\leq x_2\leq 1,\\
(\frac{C}{2}(1 - \left(\frac{\tau_y}{C}\right)^2) - \tau_y(1-\frac{\tau_y}{C}),0)^{\top}, & \textrm{ if }-\frac{\tau_y}{C}\leq x_2 \leq \frac{\tau_y}{C},\\
	(\frac{C}{2}(1-x_2^2)- \tau_y(1+x_2),0)^{\top}, & \textrm{ if }-1\leq x_2\leq-\frac{\tau_y}{C},\\
\end{array}
\right.\\
p_e(\bm{x}) := -C(x_1 - \frac{L}{2}),
\end{gather}
where $C$ is the (negative) pressure gradient. The boundary datum $\bm{u}_0$ is chosen so as to match the values in the expression above. The problem was solved with $L=4$, $C=2$ and $\tau_y=1$ using the regularization \eqref{eq:ImplicitBingham_Reg1}. Secant continuation starting from $\varepsilon = 1$ was employed to obtain better initial guesses for Newton's method; more precisely, this means that given two previously computed solutions $\bw_1,\bw_2$ corresponding to the parameters $\varepsilon_1,\varepsilon_2$, respectively, the initial guess for Newton's method at $\varepsilon$ is chosen as
\begin{equation}
	\frac{\varepsilon-\varepsilon_2}{\varepsilon_2-\varepsilon_1}(\bw_2-\bw_1) + \bw_2.
\end{equation}
In this case the tolerances were chosen to be $1\times 10^{-10}$ and $1\times 10^{-12}$ for the nonlinear and linear solvers, respectively. \textcolor{black}{Tighter tolerances are used for this problem to ensure convergence of the continuation scheme.}

\begin{figure}
\centering
\subfloat[C][{\centering Errors for different values of polynomial\par degree $k$ and number of refinements $l$.}]{{%
	\includegraphics[width=0.5\textwidth]{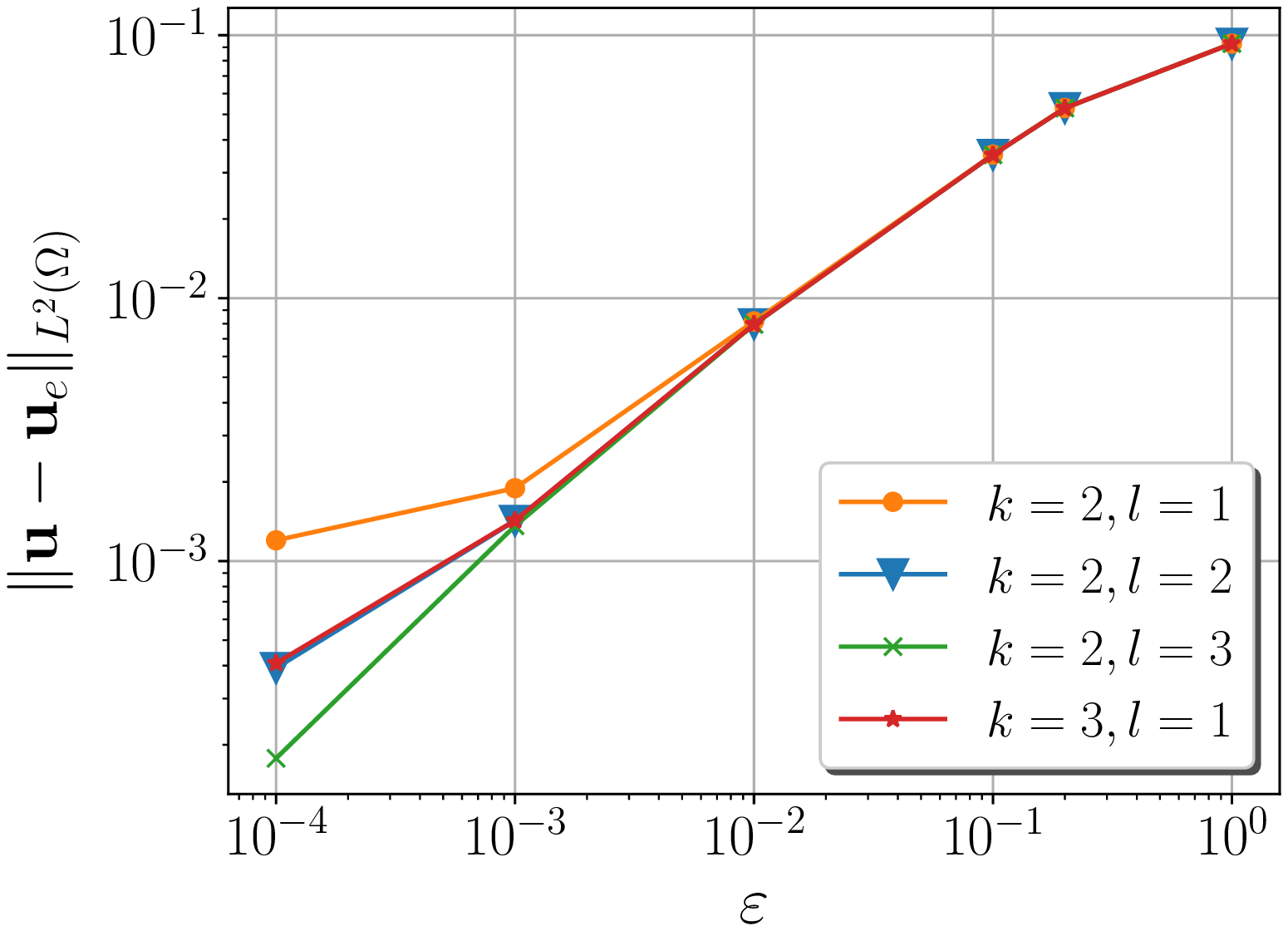}%
	}}%
	\subfloat[D][{\centering Velocity profiles at $x_1=2$ for different\par values of $\varepsilon$, including the exact solution $\varepsilon = 0$.}]{{%
	\includegraphics[width=0.5\textwidth]{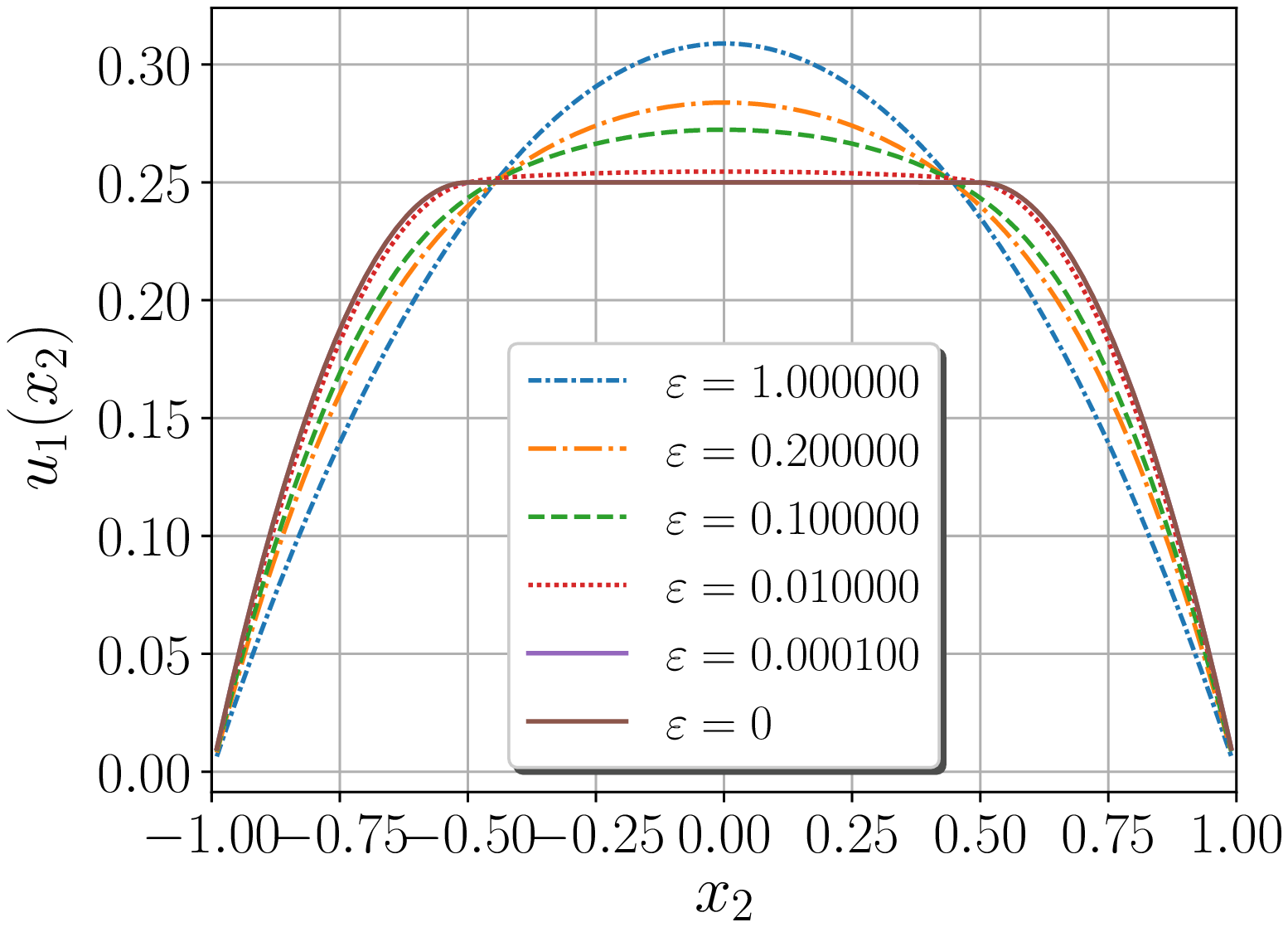}%
	}}%
	\caption{Numerical solution of the Bingham flow between two plates.}%
	\label{fig:Bingham}
\end{figure}

\cref{fig:Bingham} (a) shows the \textcolor{black}{$L^2$-distance between the numerical solution and the exact solution \cref{eq:Bingham_Poiseuille},} as $\varepsilon$ decreases, for different values of the polynomial degree $k$ and the number of refinements in the mesh hierarchy $l$; it can be observed that at some point the discretization error starts to dominate. \textcolor{black}{\cref{fig:Bingham} (b) shows the velocity profiles for different values of $\varepsilon$, including the exact solution}.  \cref{tb:BinghamPoiseuille_iterations} shows the average number of Krylov iterations per Newton step using two multigrid cycles with 5 relaxation sweeps per level as $\tilde{A}^{-1}$. It can be seen in  \cref{tb:BinghamPoiseuille_iterations} that the number of iterations remains under control, with only a slight increase for very small $\varepsilon$ and one level of refinement; the number of Newton iterations also appears to exhibit mesh-independence. In the practical computation of viscoplastic flow the approach described here should be combined with an adaptive refinement of the mesh in order to resolve the yield surface more accurately. 

\begin{table}[ht!]
\centering
\caption{\label{tb:BinghamPoiseuille_iterations}
Average number of Krylov iterations per Newton step as $\varepsilon$ decreases for the Bingham flow \protect\linebreak between two plates.}
\begin{tabular}{c c c  c c c c} 
\toprule
\multirow{2}{*}{$k$} &  \multirow{2}{*}{\# refs}& \multirow{2}{*}{\# dofs} & \multicolumn{4}{c}{$\varepsilon$ }\\ [0.1ex] 
  & & & 0.1 & 0.01 & 0.001 & 0.0001  \\
\midrule
\multirow{3}{*}{2} & 1 & $2.8\times 10^4$ & 5 & 5 & 5.33 & 14 \\
  & 2 & $1.1\times 10^5$ & 4 & 3.57 & 3.83 & 2.66 \\
  & 3 & $4.5\times 10^5$ & 4 & 4 & 3.85 & 3.5 \\ 
 \midrule
 3 & 1 & $5.9\times 10^4$ & 2.4 & 2.6 & 2.44 & 3.5 \\
  \bottomrule
\end{tabular}
\end{table}

\subsubsection*{Generalised Carreau--Yasuda Fluid}
In this experiment we employ the constitutive relation \eqref{eq:stress-Carreau} and test the solver with different values of the rheological parameters on the lid driven cavity problem. The problem is solved on the square/cube $(0,2)^d$ with $\bm{f}=\bm{0}$, and boundary data
\begin{gather*}
\bu_0(\bm{x}):= \left\{
\begin{array}{cc}
\renewcommand{\arraystretch}{1.5}
(x^2(2-x)^2,0)^\top, & \textrm{ if } y=2,\\
(0,0)^{\top}, & \textrm{ otherwise},\\
\end{array}
\right.
\end{gather*}
if $d=2$, and
\begin{gather*}
\bu_0(\bm{x}):= \left\{
\begin{array}{cc}
\renewcommand{\arraystretch}{1.5}
(x^2(2-x)^2z^2(2-z)^2,0,0)^\top, & \textrm{ if } y=2,\\
(0,0,0)^{\top}, & \textrm{ otherwise},\\
\end{array}
\right.
\end{gather*}
if $d=3$. For the 3D problem the tolerance for the linear solver was set to $1\times 10^{-8}$. In this example a simple continuation algorithm was employed to reach the different values of the parameters, e.g.\ the solution corresponding to $\nu$ is used as an initial guess in Newton's method for the problem with $\nu+\delta\nu$, iterating the procedure until the desired value is reached.  For parameters for which the effective viscosity is small (e.g.\ small $\nu$), the problem will be convection dominated and hence some advective stabilization is required in \eqref{eq:FEM_nonlinear_momentum}. We choose to add a stabilizing term based on jump penalisation described in \cite{Burman2008,Douglas1976}:
\begin{equation}
	S_h(\bv,\bw) := \sum_{K\in\mathcal{M}_h} \frac{1}{2}\int_{\partial K}\delta \, h^2_{\partial K}\,\jump{\nabla\bv} : \jump{\nabla\bw},
\end{equation}
where $\jump{\bz}$ denotes the jump of $\bz$ across $\partial K$,  $h_{\partial K}$ is a function giving the size of each face in $\partial K$, and $\delta$ is an arbitrary stabilization parameter. In the numerical experiments the stabilization parameter was chosen to be cell-dependent and set to $5\times10^{-3}\|\tilde{u}\|_{L^\infty(K)}$. In the experiments described in this section, 2 full multigrid cycles with 4 relaxation sweeps per level were applied as $\tilde{A}^{-1}$ when $d=2$, and 1 cycle with 6 relaxation sweeps when $d=3$\textcolor{black}{. These values were chosen so as to balance the amount of inner and outer work (e.g.\ fewer relaxation sweeps result in less expensive linear solves, but more iterations are needed); convergence is also achieved with fewer relaxation sweeps, but the values chosen here resulted in a shorter time to solution}. \cref{tb:Carreau_iterations_nu,tb:Carreau_iterations_nu_3D} show the average number of Krylov iterations per Newton step for a problem with decreasing $\nu$; it can be observed that the number of iterations remains well controlled even for the lowest values of $\nu$ (in the Newtonian problem, $\nu=0.0002$ would correspond to a Reynolds number of $10000$).

\begin{table}[ht!]
\centering
\caption{\label{tb:Carreau_iterations_nu}
Average number of Krylov iterations per Newton step as $\nu$ decreases for the 2D generalised Carreau--Yasuda relation with $r_1=1.8$, $r_2 =2.5$, $\Gamma_1=\Gamma_2=200$, $\beta_1=0.9$, $\beta_2 = 0.5$.}
\begin{tabular}{c c c  c c c c} 
\toprule
\multirow{2}{*}{$k$} &  \multirow{2}{*}{\# refs}& \multirow{2}{*}{\# dofs} & \multicolumn{4}{c}{$\nu$ }\\ [0.1ex] 
  &  & & 0.2 & 0.001 & 0.0005 & 0.0002  \\
\midrule
\multirow{3}{*}{2} & 1 & $3.1\times 10^4$ & 4.25 & 3.5 & 4 & 5 \\
  & 2 & $1.2\times 10^5$ & 4.25 & 3.5 & 3.5 & 4 \\
  & 3 & $4.9\times 10^5$ & 4.25 & 3 & 2.5 & 3 \\ 
 \midrule
\multirow{3}{*}{3} & 1 & $6.5\times 10^4$ & 2.75 & 2. & 2.5 & 2.5 \\
  & 2 & $2.5\times 10^5$ & 2.75 & 1.66 & 2 & 2.5 \\
  & 3 & $1.0\times 10^6$ & 2.5 & 2 & 1.5 & 1.5 \\ 
  \bottomrule
\end{tabular}\\
\end{table}

\begin{table}[ht!]
\centering
\caption{\label{tb:Carreau_iterations_nu_3D}
Average number of Krylov iterations per Newton step as $\nu$ decreases for the 3D generalised Carreau--Yasuda relation with $r_1=1.8$, $r_2 =2.5$, $\Gamma_1=\Gamma_2=200$, $\beta_1=0.9$, $\beta_2 = 0.5$.}
\begin{tabular}{c c c  c c c c} 
\toprule
\multirow{2}{*}{$k$} &  \multirow{2}{*}{\# refs}& \multirow{2}{*}{\# dofs} & \multicolumn{4}{c}{$\nu$ }\\ [0.1ex] 
  &  & & 0.2 & 0.002 & 0.0005 & 0.00028  \\
\midrule
3 & 1 & $9.2\times 10^5$ & 7.25 & 5 & 5.5 & 5.5 \\
  \bottomrule
\end{tabular}\\
\end{table}

\textcolor{black}{A comparison with the preconditioner using a Jacobi smoother instead of the macrostar iteration can be found in \cref{tb:timing}, for a given set of rheological parameters. The experiments were performed on 12 Intel Xeon Silver 4116 CPUs. Very mild parameters are considered for this comparison, since the Krylov solver with Jacobi smoothing fails to converge otherwise (the solver using the AMG libraries Hypre \cite{falgout2002hypre}, ML \cite{Gee2006}, and GAMG \cite{adams2004} on the stress-velocity block failed to converge altogether). We note that for our academic test problems, the preconditioner employing a direct sparse solver for the stress-velocity block is still faster on the workstation resources we had available, but we expect that the implementation could be optimised and the algorithm employing the macrostar iteration will scale better on high performance computers. Other ways of lowering the cost of the algorithm, such as employing $H(\diver)$--$L^2$-type elements for the velocity and pressure, for which a smaller star iteration would suffice to capture the kernel, will be the subject of future research.}


\begin{table}[ht!]
\centering
\caption{\label{tb:timing}
Runtime comparison and total number of Krylov iterations (\# iters) for a 2D problem with $\nu=0.02$, $r_1=1.8$, $r_2=2.5$, $\Gamma_1=\Gamma_2=200$, $\beta_1=0.9$, $\beta_2=0.5$, and a 3D problem with $r_1=r_2=\nu =2$. }
\begin{tabular}{c c c  c c c c} 
\toprule
\multirow{2}{*}{$d$} &  \multirow{2}{*}{$(k,\mathrm{\# refs})$}& \multirow{2}{*}{\# dofs} & \multicolumn{2}{c}{macrostar} & \multicolumn{2}{c}{Jacobi} \\ [0.1ex] 
		     &  & & \# iters & time (min.) & \# iters & time (min.)  \\
\midrule
2 & $(2,3)$ & $4.9\times 10^5$ & 15 & 1.67 & 3040 & 107.81 \\
3 & $(3,1)$ & $9.2\times 10^5$ & 10 & 37.33 & 753 & 70.35 \\
  \bottomrule
\end{tabular}\\
\end{table}

\cref{tb:Carreau_iterations_r,tb:Carreau_iterations_eps} show the number of average Krylov iterations for small $r_1$ and large $\Gamma_2$, respectively, for two different values of $\gamma$. It can be observed that depending on the parameter of interest, large values of $\gamma$ improve the robustness of the algorithm.
In all the examples \textcolor{black}{in this section,} the solver appears to be robust with respect to the parameters appearing in the constitutive relation and also exhibits mesh-independence.

\begin{table}[ht!]
\centering
\caption{\label{tb:Carreau_iterations_r}
Average number of Krylov iterations per Newton step as $r_1$ decreases for the 2D generalised Carreau--Yasuda relation with $\nu = 0.01$, $r_2 =2$, $\Gamma_1=125$, $\beta_1=0.7$.}
\begin{tabular}{c c c c  c c c c} 
\toprule
 \multirow{2}{*}{$\gamma$} & \multirow{2}{*}{$k$} &  \multirow{2}{*}{\# refs}& \multirow{2}{*}{\# dofs} & \multicolumn{4}{c}{$r_1$ }\\ [0.1ex] 
& &  & & 1.66 & 1.25 & 1.11 & 1.07  \\
\midrule
\multirow{6}{*}{$10^4$}& \multirow{3}{*}{2} & 1 & $3.1\times 10^4$ & 3.5 & 3.5 & 3.5 & 3.5 \\
&  & 2 & $1.2\times 10^5$ & 3.5 & 3.5 & 3.5 & 3.5 \\
&  & 3 & $4.9\times 10^5$ & 3 & 3.5 & 4 & 4 \\ 
\cmidrule{3-8}
& \multirow{3}{*}{3} & 1 & $6.5\times 10^4$ & 2 & 2 & 2 & 2 \\
&  & 2 & $2.5\times 10^5$ & 2 & 2 & 2 & 2.5 \\
&  & 3 & $1.0\times 10^6$ & 2 & 2 & 2.5 & 2.5 \\ 
\midrule
\multirow{6}{*}{$1$}& \multirow{3}{*}{2} & 1 & $3.1\times 10^4$ & 5 & 4 & 4 & 4 \\
&  & 2 & $1.2\times 10^5$ & 4.5 & 4 & 3.5 & 3.5 \\
&  & 3 & $4.9\times 10^5$ & 4 & 4 & 4 & 4 \\ 
\cmidrule{3-8}
& \multirow{3}{*}{3} & 1 & $6.5\times 10^4$ & 4 & 4 & 3.5 & 3 \\
&  & 2 & $2.5\times 10^5$ & 4 & 3.5 & 3 & 3 \\
&  & 3 & $1.0\times 10^6$ & 4 & 3.5 & 3 & 3 \\ 
  \bottomrule
\end{tabular}\\
\end{table}
 \begin{remark}
 In general, extreme values of the parameters could result in convergence issues for the nonlinear iterations. In practice, the preconditioner presented here should then be coupled e.g.\ with a more sophisticated continuation strategy for the nonlinear iterations, or with nested iteration. 
 \end{remark}

\begin{table}[ht!]
\centering
\caption{\label{tb:Carreau_iterations_eps}
Average number of Krylov iterations per Newton step as $\Gamma_2$ increases for the 2D generalised Carreau--Yasuda relation with $\nu=0.01$, $r_1=1.7$, $r_2 =3$, $\Gamma_1=10$, $\beta_1=0.2$, $\beta_2 = 0.9$.}
\begin{tabular}{c c c c  c c c c} 
\toprule
 \multirow{2}{*}{$\gamma$} & \multirow{2}{*}{$k$} &  \multirow{2}{*}{\# refs}& \multirow{2}{*}{\# dofs} & \multicolumn{4}{c}{$\Gamma_2$ }\\ [0.1ex] 
& &  & & 10 & 1000 & 5000 & 10000  \\
\midrule
\multirow{6}{*}{$10^4$} & \multirow{3}{*}{2} & 1 & $3.1\times 10^4$  & 4 & 3.66 & 4 & 4 \\
&  & 2 & $1.2\times 10^5$  & 3.66 & 3.66 & 4 & 3.5 \\
&  & 3 & $4.9\times 10^5$ & 3.66 & 3.66 & 4 & 4 \\ 
\cmidrule{3-8}
& \multirow{3}{*}{3} & 1 & $6.5\times 10^4$ & 2.33 & 2 & 2.5 & 2.5 \\
&  & 2 & $2.5\times 10^5$ & 2.33 & 2 & 2.5 & 2.5 \\
&  & 3 & $1.0\times 10^6$ & 2.33 & 2 & 2.5 & 2.5 \\ 
\midrule
\multirow{6}{*}{$1$}& \multirow{3}{*}{2} & 1 & $3.1\times 10^4$ & 9.66 & 20.3 & 32 & 34.5 \\
&  & 2 & $1.2\times 10^5$ & 9 & 19.3 & 30.5 & 31.5 \\
&  & 3 & $4.9\times 10^5$ & 8 & 17.3 & 26 & 27 \\ 
\cmidrule{3-8}
& \multirow{3}{*}{3} & 1 & $6.5\times 10^4$ & 7.33 & 15.6 & 24 & 26 \\
&  & 2 & $2.5\times 10^5$ & 6.33 & 13.6 & 20 & 21 \\
&  & 3 & $1.0\times 10^6$ & 6 & 11.3 & 16.5 & 17.5 \\ 
  \bottomrule
\end{tabular}\\
\end{table}

\subsubsection*{Activated Euler Flow Past an Obstacle}
Consider the non-standard constitutive relation 
\begin{equation}\label{eq:ActivatedFluid}
\renewcommand{\arraystretch}{1.5}
\left\{
\begin{array}{cc}
	\BD = \tau_y\frac{\BS}{|\BS|} + \BD_2, & \textrm{ if } |\BD|\geq \tau_y,\\
  \BS = 0, & \textrm{ if }|\BD|<\tau_y,\\
\end{array}
\right.
\end{equation}
where $\BD_2$ satisfies $\BS = 2\nu |\BD_2|^{r-2}\BD_2$, for some $\nu>0$ and $r>1$. A material with a response of this type of response is called an \emph{Euler/power-law fluid} and some of its properties were analyzed for the first time in \cite{Blechta2019}; it describes an inviscid fluid before activation (i.e.\ when $|\BD|\leq \tau_y$) and a power-law fluid otherwise. It is likely that such constitutive relations had not been considered before due to the prevalence of explicit relations of the kind $\BS = \BS(\BD)$ in the literature, but they could nevertheless be potentially useful in applications and further research is warranted.

\begin{figure}
\centering
\subfloat[E][{Effective viscosity for an activated Euler/power-law fluid with $r=1.3$, $\tau_y = 3$, $\nu = 0.5$, and $\varepsilon= 1\times 10^{-5}$.}]{{%
	\includegraphics[width=\textwidth]{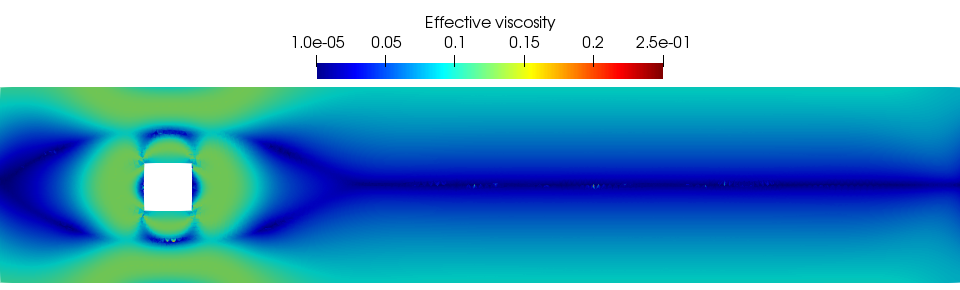}%
}}\\
\subfloat[F][{Effective viscosity for a power-law fluid with $r=1.3$ and $\nu = 0.5$.}]{{%
	\includegraphics[width=\textwidth]{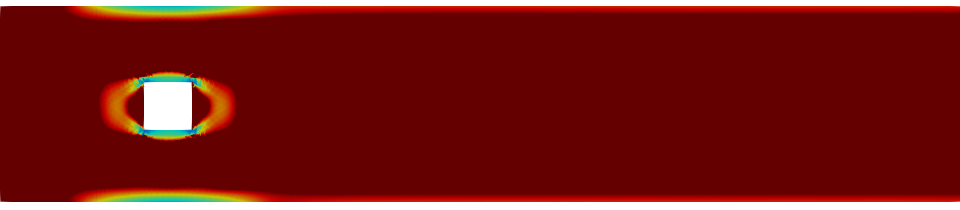}%
}}%
	\caption{Effective viscosity for the flow past an obstacle.}%
	\label{fig:ActivatedEuler}%
\end{figure}

Observe that the power-law nonlinearity can be inverted and we have that, for any $\BD,\BS\in\Rds$,
\begin{equation}
	\BS = 2\nu |\BD|^{r-2}\BD  \Longleftrightarrow  \BD = \frac{1}{2\nu} \left(\frac{|\BS|}{2\nu} \right)^{r'-2}\BS.
\end{equation}
Using this fact we can write a regularized constitutive relation similar to the one described in  \cref{rem:ActEuler}:
\begin{equation*}
	\BG(\BS,\BD) = \BD - \left(\frac{1}{2\nu}\left|\frac{\BS}{2\nu}\right|^{r'-2} + \frac{\tau_y}{\sqrt{\varepsilon^2 + |\BS|^2}} \right)\BS,
\end{equation*}
where $\varepsilon>0$. The problem was solved on the set $\Omega = (0,2)\times(0,0.41) \setminus (0.3,0.4)\times (0.15,0.25)$, with boundary data
\begin{equation}
\renewcommand{\arraystretch}{1.5}
\left\{
\begin{array}{cc}
	\bu  = (4\frac{0.3 x_2(0.41-x_2)}{0.41^2},0)^\top, & \textrm{ on } \partial\Omega_1 := \{x_1 = 0\}\cap\partial\Omega,\\
	\bu_\tau = 0 \textrm{ and } \BS\bm{n}\cdot\bm{n} - p = 0 & \textrm{ on }\partial\Omega_2 := \{x_1 = 2\}\cap\partial\Omega,\\
	\bu  = \bm{0}, & \textrm{ on } \partial\Omega\setminus(\partial\Omega_1\cup\partial\Omega_2),\\
\end{array}
\right.
\end{equation}
where $\bm{n}$ is the outward normal vector to the boundary and $\bm{u}_\tau = \bu - (\bu\cdot\bm{n})\bm{n}$ is the tangential part of the velocity. \cref{tb:Activated_thin} shows the number of Krylov iterations per Newton step obtained using two full multigrid cycles with 3 relaxation steps per level as $\tilde{A}^{-1}$; the same robust behaviour as in the previous examples can be observed here.  \cref{fig:ActivatedEuler} shows the effective viscosity $\mu_{\textrm{eff}} := \frac{\BS}{2\BD}$ for the solution of this problem and for that of a regular shear-thinning power-law fluid. It can be observed that the effective viscosity of the activated fluid greatly decreases far away from the obstacle, which is a common assumption in the study of boundary layers.

\begin{table}[ht!]
\centering
\caption{\label{tb:Activated_thin}
Average number of Krylov iterations per Newton step as $\varepsilon$ decreases for the Euler/power-law relation with $\nu=0.5$, $r=1.3$, $\tau_y=3$.}
\begin{tabular}{c c c  c c c c} 
\toprule
\multirow{2}{*}{$k$} &  \multirow{2}{*}{\# refs}& \multirow{2}{*}{\# dofs} & \multicolumn{4}{c}{$\varepsilon$ }\\ [0.1ex] 
  & & & 0.2 & 0.01 & 0.0001 & 0.00001  \\
\midrule
\multirow{3}{*}{2} & 1 & $3.5\times 10^4$ & 5 & 3 & 2 & 2 \\
  & 2 & $1.4\times 10^5$ & 5.66 & 4 & 2 & 2 \\
  & 3 & $5.6\times 10^5$ & 4.6 & 4 & 3 & 3 \\ 
 \midrule
\multirow{3}{*}{3} & 1 & $7.3\times 10^4$ & 2.66 & 2 & 1 & 1 \\
  & 2 & $2.9\times 10^5$ & 3 & 2 & 2 & 2 \\
  & 3 & $1.1\times 10^6$ & 3 & 2 & 2 & 2 \\ 
  \bottomrule
\end{tabular}\\
\end{table}

\section{Conclusion}
In this work we have extended the work of \cite{Benzi2006,Farrell2019a,Farrella,Farrell2020}
on parame\-ter-robust preconditioners for two-field formulations of the
Navier--Stokes equations to a three-field formulation of the non-Newtonian equations
with an implicit constitutive relation. An augmented Lagrangian term controls
the Schur complement with respect to the pressure, while a specialized multigrid
scheme is applied monolithically to the augmented stress-velocity block. The
preconditioner is robust to variation of the rheological parameters in numerical
experiments. We expect that the same strategy will apply straightforwardly to
the transient case \cite{Heister2013}. An important extension is to the
anisothermal case: future work will investigate whether it is advantageous to
take the Schur complement with respect to the temperature \cite{Howle2012}, or
treat the temperature monolithically with the stress and velocity.

\bibliographystyle{siamplain}
\bibliography{bibliography}
\end{document}